\documentclass[final]{siamltex}

\usepackage{latexsym,bm}
\usepackage{mathrsfs,amsmath,amssymb}
\usepackage{amscd} 
\usepackage{relsize,amsmath} 

\newcommand{\lint}{\mathlarger{\int}}
\newcommand{\Lint}{\mathlarger{\mathlarger{\int}}}

\makeatletter
\newcommand{\biggg}[1]{{\hbox{$\left#1\vbox to 20.5pt{}\right.\n@space$}}}
\newcommand{\Biggg}[1]{{\hbox{$\left#1\vbox to 23.5pt{}\right.\n@space$}}}
\newcommand{\bigggg}[1]{{\hbox{$\left#1\vbox to 26.5pt{}\right.\n@space$}}}
\newcommand{\Bigggg}[1]{{\hbox{$\left#1\vbox to 29.5pt{}\right.\n@space$}}}
\newcommand{\biggggg}[1]{{\hbox{$\left#1\vbox to 32.5pt{}\right.\n@space$}}}
\newcommand{\Biggggg}[1]{{\hbox{$\left#1\vbox to 35.5pt{}\right.\n@space$}}}
\newcommand{\bigggggg}[1]{{\hbox{$\left#1\vbox to 38.5pt{}\right.\n@space$}}}
\newcommand{\Bigggggg}[1]{{\hbox{$\left#1\vbox to 41.5pt{}\right.\n@space$}}}
\makeatother
\newcommand{\bigggl}{\mathopen\biggg}

\newcommand{\bigggr}{\mathclose\biggg}
\newcommand{\Bigggl}{\mathopen\Biggg}

\newcommand{\Bigggr}{\mathclose\Biggg}

\DeclareMathOperator{\wlim}{w-\lim}
\DeclareMathOperator{\w*lim}{w^{*}-\lim}
\DeclareMathOperator{\slim}{s-\lim}

\title{Asymptotic behavior of solutions toward \\
a multiwave pattern 
to the Cauchy problem \\
for the scalar conservation law \\
with degenerate flux and viscosity}

\author{Natsumi Yoshida \thanks{
BKC Research Organization of Social Sciences, 
Ritsumeikan University, Kusatsu, Shiga 525-8577, Japan
({\tt 14v00067@gst.ritsumei.ac.jp})/
Osaka City University Advanced Mathematical Institute, 
Sumiyoshi, Osaka 558-8585, Japan.}}

\begin{document}

\maketitle

\begin{abstract}
In this paper, we investigate the asymptotic behavior of solutions toward 
a multiwave pattern of the Cauchy problem 
for the scalar viscous conservation law 
where the far field states are prescribed. 
Especially, we deal with the case 
when the flux function is convex or concave 
but linearly degenerate on some interval, 
and also the viscosity 
is a nonlinearly degenerate one ($p$-Laplacian type viscosity). 
When the corresponding Riemann problem admits a Riemann solution which 
consists of rarefaction waves and contact discontinuity, 
it is proved that 
the solution of the Cauchy problem tends toward the linear combination of 
the rarefaction waves and contact wave for $p$-Laplacian type viscosity 
as the time goes to infinity. 
This is the first result concerning the asymptotics 
toward multiwave pattern 
for the Cauchy problem of the scalar conservation law 
with nonlinear viscosity. 
The proof is given by a technical energy methods 
and the careful estimates 
for the interactions between the nonlinear waves. 
\end{abstract}



\pagestyle{myheadings}
\thispagestyle{plain}
\markboth{N. YOSHIDA}{STABILITY OF MULTIPLE NONLINEAR WAVES}

\section{Introduction and main theorem}
In this paper, 
we shall consider the asymptotic behavior of 
solutions 
for one-dimensional scalar conservation law 
with a nonlinearly degenerate viscosity 
($p$-Laplacian type viscosity with $p>1$) 
\begin{eqnarray}
 \left\{\begin{array}{ll}
  \partial_tu +\partial_x \bigl(f(u) \bigr)
  = \mu \, 
    \partial_x \left( \, 
    \left| \, \partial_xu \, \right|^{p-1} \partial_xu \, 
    \right)
  \qquad &(t>0, x\in \mathbb{R}), \\[5pt]
  u(0,x) = u_0(x) \qquad &( x \in \mathbb{R} ),\\[5pt]
  \displaystyle{\lim_{x\to \pm \infty}} u(t,x) =u_{\pm}  
  \qquad &\bigl( t \ge 0 \bigr).   
 \end{array}
 \right.\,
\end{eqnarray}
Here, $u=u(t,x)$ denotes the unknown function of $t>0$ and $x\in \mathbb{R}$, 
the so-called conserved quantity, 
$f=f(u)$ is the flux function depending only on $u$, 
$\mu$ is the viscosity coefficient, 
$u_0$ is the given initial data, 
and constants $u_{\pm } \in \mathbb{R}$ 
are the prescribed far field states. 
We suppose the given flux $f=f(u)$ is a $C^1$-function 
satisfying $f(0)=f'(0)=0$, 
$\mu$ is a positive constant 
and far field states $u_{\pm }$ satisfy $u_-<u_+$ 
without loss of generality. 

We are interested in the asymptotic behavior 
and its precise estimates 
in time of 
the global solution to the Cauchy problem (1.1). 
Especially, 
one of the keys of the study is to investigate 
the influence of 
the shape of the flux function $f(u)$ and the far field states $u_{\pm }$ 
on the asymptotic behavior. 
It can be expected that the large-time behavior 
is closely related to the weak solution (``Riemann solution'') of the corresponding 
Riemann problem 
(cf. \cite{lax}, \cite{smoller})
for the non-viscous hyperbolic part of (1.1): 
\begin{eqnarray}
 \left\{\begin{array} {ll}
 \partial _t u + \partial _x \bigl( f(u) \bigr)=0 
 \qquad &(t>0, x\in \mathbb{R}),\\[5pt]
u(0,x)=u_0 ^{\rm{R}} (x)\qquad &(x \in \mathbb{R}),
 \end{array}
  \right.\,
\end{eqnarray}
where $u_0 ^{\rm{R}}$ is the Riemann data defined by
$$
u_0 ^{\rm{R}} (x)=u_0 ^{\rm{R}} (x\: ;\: u_- ,u_+)
          := 
          \left\{\begin{array} {ll}
          u_-  & \; (x < 0),\\[5pt]
          u_+  & \; (x > 0).
          \end{array}\right.
$$
In fact, for the usual linear viscosity 
case:
\begin{eqnarray}
 \left\{\begin{array}{ll}
  \partial_tu +\partial_x \bigl(f(u) \bigr)= \mu \, \partial_x^2 u
  \qquad &(t>0, x\in \mathbb{R}), \\[5pt]
  u(0,x) = u_0(x) \qquad &( x \in \mathbb{R} ),\\[5pt]
  \displaystyle{\lim_{x\to \pm \infty}} u(t,x) =u_{\pm}  
  \qquad &\bigl( t \ge 0 \bigr),   
 \end{array}
 \right.\,
\end{eqnarray}
when the smooth flux function $f$ is genuinely nonlinear 
on the whole space $\mathbb{R}$,
i.e., $f''(u)\ne 0\ (u\in \mathbb{R})$, 
Il'in-Ole{\u\i}nik \cite{ilin-oleinik} showed 
the following: if  $f''(u) > 0\ (u\in \mathbb{R})$, that is, the Riemann solution 
consists of a single rarefaction wave solution, 
the global solution in time of the Cauchy problem 
(1.3) tends toward the rarefaction wave;
if  $f''(u) < 0\ (u\in \mathbb{R})$, that is,
the Riemann solution consists of a single shock wave solution, 
the global solution of the
Cauchy problem (1.3) does the corresponding smooth traveling wave solution 
(``viscous shock wave'') of (1.3) with a spacial shift 
(cf. \cite{ilin-kalashnikov-oleinik}). 
More generally, in the case of the flux functions 
which are not uniformly genuinely nonlinear,
when the Riemann solution consists of a single shock wave 
satisfying 
Ole{\u\i}nik's shock condition, 
Matsumura-Nishihara \cite{matsu-nishi3} showed the asymptotic stability 
of the corresponding viscous shock wave. 
However, 
when we consider 
the circumstances 
where the Riemann solution generically forms a pattern of
multiple nonlinear waves which consists of 
rarefaction waves, shock waves 
and waves of contact discontinuity (refer to \cite{liep-rosh}), 
there had been no results about the asymptotics toward 
the multiwave pattern. 
Recently, Matsumura-Yoshida \cite{matsumura-yoshida} 
proved the asymptotics toward 
a multiwave pattern 
of the superposition of the rarefaction waves 
and a self-similar solution (``viscous contact wave'') 
which is corresponded to the wave of the contact discontinuity. 
Namely, they investigated 
the case where the flux function $f$ is smooth
and genuinely nonlinear (that is, $f$ is convex function or concave function) 
on the whole $\mathbb{R}$ except 
a finite interval $I := (a,b) \subset \mathbb{R}$, and 
linearly degenerate on $I$, that is, 
\begin{equation}
\left\{
\begin{array}{ll}
  f''(u) >0 & \; \bigl(u \in (-\infty ,a\, ]\cup [\, b,+\infty )\bigr),\\[5pt]
  f''(u) =0 & \; \bigl(u \in (a,b)\bigr).
\end{array}\right.
\end{equation}
For 
the flux function satisfying (1.4),  
the corresponding Riemann solution does form 
multiwave pattern which consists of the contact discontinuity 
with the jump from $u=a$ to $u=b$ and the rarefaction waves, 
depending on the choice of $a$, $b$, $u_-$ and $u_+$.
Thanks to that the cases in which the interval $(a,b)$ is disjoint 
from the interval $(u_-,u_+)$ 
are 
similar as in the case 
the flux function $f$ is genuinely nonlinear 
on the whole space $\mathbb{R}$, 
and the case $u_-<a<u_+<b$ is 
the same as that for $a<u_-<b<u_+$, 
we may only consider the 
typical cases 
\begin{equation}
a<u_-<b<u_+ \quad \mbox{or} \quad u_-<a<b<u_+.   
\end{equation}
Under the conditions (1.4) and (1.5), 
they have shown the unique global solution in time 
to the Cauchy problem (1.3) tends uniformly in space toward 
the multiwave pattern of the combination of 
the viscous contact wave 
and the rarefaction waves 
as the time goes to infinity. 
It should be noted that the rarefaction wave which connects the 
far field states $u_-$ and $u_+$ 
$\bigl(u_\pm \in (-\infty,a\, ]$ or $u_\pm \in [\, b,\infty)\bigr)$
is explicitly given by
\begin{equation}
u=u^r \left( \frac{x}{t}\: ;\:  u_- , u_+ \right)
:= \left\{
\begin{array}{ll}
  u_-  & \; \bigl(\, x \leq \lambda(u_-)\,t \, \bigr),\\[7pt]
  \displaystyle{ (\lambda)^{-1}\left( \frac{x}{t}\right) } 
  & \; \bigl(\, \lambda(u_-)\,t \leq x \leq \lambda(u_+)\,t\,  \bigr),\\[7pt]
   u_+ & \; \bigl(\, x \geq \lambda(u_+)\,t \, \bigr),
\end{array}
\right.
\end{equation} 
where $\lambda(u):=f'(u)$,
and the viscous contact wave 
which connects 
$u_-$ and $u_+$ ($u_\pm \in [\,a,b\,]$) is given by an exact
solution of the linear 
convective heat equation 
\begin{equation}
  \partial _t u + \tilde{\lambda} \, \partial_x u = \mu \, \partial_x^2 u
  \qquad \biggl(
  \,\tilde{\lambda}:= \frac{f(b)-f(a)}{b-a},\ 
  t>0,x \in  \mathbb{R} 
  \biggr)
\end{equation}
 which has the form
$$
u=U\left(\frac{x-\tilde{\lambda} \, t}{\sqrt{t}}\: ;\: u_- ,u_+ \right)
$$
where $U\left(\frac{x}{\sqrt{t}}\,;\,u_- ,u_+ \right)$ is explicitly defined by
\begin{equation}
 U\left(\frac{x}{\sqrt{t}}\: ;\: u_- ,u_+ \right)
 :=u_- +\frac{u_+ - u_-}{\sqrt{\pi}}
   \lint ^{\frac{\mathlarger{x}}{\mathlarger{\sqrt{4\mu t}}}}_{-\infty} 
 \mathrm{e}^{-\xi^2}\, \mathrm{d}\xi 
\quad \, (t>0, x\in \mathbb{R}). 
\end{equation}
Yoshida \cite{yoshida} also obtained the precise decay properties 
for the asymptotics (cf. \cite{hashimoto-kawashima-ueda}). 
In the proof of them, the a priori energy estimates 
acquired by an $L^2$-energy method and careful 
estimates for the terms of nonlinear interactions of 
the viscous contact wave and the rarefaction waves. 

The aim of the present paper is 
to extend the results in the previous study in \cite{matsumura-yoshida} 
to the case where the viscosity is of $p$-Laplacian 
type 
(the related problems are studied in 
\cite{gur-mac}, \cite{nagai-mimura1}, \cite{nagai-mimura2} 
and so on). 
For this case, a main difficulty arises from the fact that 
when $u_\pm \in [\, a,b\, ]$, 
the asymptotic state is expected to be a self-similar type solution 
of a nonlinearly degenerate convective heat equation 
which may need the more subtle treatment than 
the Gaussian kernel type one (1.8) of the equation (1.7). 
There is only one result for the asymptotic behavior for the problem (1.1) 
in the case where the flux function is genuinely nonlinear 
on the whole space $\mathbb{R}$. 
Namely, Matsumura-Nishihara \cite{matsu-nishi2} 
proved the asymptotics which tends toward a single rarefaction wave 
by using the $L^2$ and $L^p$-energy estimates. 
%
We then consider the case where the flux function is given as (1.4) 
and the far field states as (1.5). 
We expect the asymptotic behavior of solutions to the Cauchy problem (1.1) 
to be similar as in \cite{matsumura-yoshida}. 
In more detail, under the conditions (1.4) and (1.5), 
if the far field states $u_{\pm}$ satisfy 
$u_\pm \in (-\infty,a\, ]$ or $u_\pm \in [\, b,\infty)$, 
the asymptotic state of the solutions to the Cauchy problem (1.1) should be the rarefaction wave (1.6) 
which connects $u_-$ and $u_+$, 
and if the far field states $u_{\pm}$ satisfy 
$u_\pm \in [\, a,b\, ]$
, 
the one should be the ``contact wave for $p$-Laplacian type viscosity'' 
which connects $u_-$ and $u_+$, which is given by an exact solution of the following 
$p$-Laplacian evolution equation 
\begin{equation}
  \partial _t u + \tilde{\lambda} \, \partial_x u 
  = \mu \, 
    \partial_x \left( \, 
    \left| \, \partial_xu \, \right|^{p-1} \partial_xu \, 
    \right)
  \qquad \biggl(
  \,\tilde{\lambda}:= \frac{f(b)-f(a)}{b-a},\ 
  t>0,x \in  \mathbb{R} 
  \biggr). 
\end{equation}
In order to look for an exact solution, 
especially self-similar type solution, 
we differentiate the evolution equation (1.9) 
with respect to $x$ 
and we have the following porous medium equation 
with the convection term
\begin{equation}
\partial_t v  + \tilde{\lambda} \, \partial_x v 
= \mu \, 
  \partial_x^2 
  \left( \, 
  \left| \, v \, \right|^{p-1} v \, 
  \right), 
\end{equation}
where $v:=\partial_x u$. 
Barenblatt \cite{barenblatt}, 
Zel'dovi{\v{c}}-Kompanceec \cite{zel-kom} and Pattle \cite{pattle} 
(see also \cite{carillo-toscani}, \cite{huang-pan-wang}, \cite{kamin}) 
introduced the following Cauchy problem of the porous medium equation 
\begin{eqnarray}
 \left\{\begin{array}{ll}
  \partial_t v
  = 
    \mu \, 
    \partial_x^2 
    \left( \, 
    \left| \, v \, \right|^{p-1} v \, 
    \right)
  \qquad &(t>-1, x\in \mathbb{R}), \\[5pt]
  v(-1,x) = (u_{+} - u_{-})\, \delta(x) 
  \qquad &( x \in \mathbb{R}\, ;\,  u_{-}<u_{+}),\\[5pt]
  \displaystyle{\lim_{x\to \pm \infty}} v(t,x) =0  
  \qquad &\bigl( t \ge -1 \bigr),  
 \end{array}
 \right.\,
\end{eqnarray}
where $\delta(x)$ is the Dirac $\delta$-distribution. 
They obtained the Barenblatt-Kompanceec-Zel'dovi{\v{c}} solution 
\begin{equation}
v(t,x)
:= \frac{1}{(1+t)^{\frac{1}{p+1}}}\, 
   \Bigggl( \left( 
   A-B \left(\frac{x}{(1+t)^{\frac{1}{p+1}}} \right)^2 \, 
   \right)\vee 0 \Bigggr)^{\frac{1}{p-1}},
\end{equation}
\begin{eqnarray*}
\left\{\begin{array}{ll}
\displaystyle{
A=A_{p,\mu,u_{\pm}}
:= \left( \, 
   \frac{(p-1)\, \left( \, u_{+} - u_{-} \, \right)}
   {8\, \mu\, p(p+1) 
   \biggl( \, \displaystyle{\int^{\frac{\pi}{2}}_{0} 
   \bigl( \sin \theta \bigr)^{\frac{p+1}{p-1}}\, \mathrm{d}\theta }
   \biggr)^2 }
   \, \right)^{\frac{p-1}{p+1}}
}, \\[25pt]
\displaystyle{ 
B=B_{p,\mu}:= \frac{p-1}{2\, \mu \, p(p+1)} 
}, \\[10pt] 
\displaystyle{ 
2A^{\frac{p+1}{2(p-1)}}B^{-\frac{1}{2}}
      \int^{\frac{\pi}{2}}_{0} 
      \bigl( \sin \theta \bigr)^{\frac{p+1}{p-1}}\, \mathrm{d}\theta 
      = u_{+} - u_{-}},
\end{array}
\right.\,
\end{eqnarray*}
where the symbol ``$\vee $'' is defined as 
$
a\vee b:= \max \{a,b\}. 
$
Thus, when we define by using the solution (1.12) as 
\begin{align}
\begin{aligned}
U\left(\, \frac{x}{(1+t)^{\frac{1}{p+1}}}\: ;\: u_- ,u_+ \right)
&:= u_- + 
    \int^{x}_{-\infty} v(t,y)\, \mathrm{d}y \\
&=u_- + \displaystyle{
        \lint^{\frac{\mathlarger{x}}{ \mathlarger{(1+t)}^{{\scriptscriptstyle \frac{1}{p+1}}} }}_{-\infty}
             \Bigl( \left( \, 
                A-B \xi^2 \, 
                \right)\vee 0 \, \Bigr)^{\frac{1}{p-1}}
                \, \mathrm{d}\xi } 
             \\
&\quad \, (t>-1, x\in \mathbb{R}), 
\end{aligned}
\end{align}
and change the variable as $1+t \mapsto t > 0$
and $x \mapsto x-\tilde{\lambda} \, t$ in this order, 
we have a desired 
canditate of the asymptotic state as 
\begin{align}
U\left(\, \frac{x-\tilde{\lambda} \, t}
          {t^{\frac{1}{p+1}}}\: ;\: u_- ,u_+ \right) 
= u_- + \displaystyle{
        \lint^{\frac{\mathlarger{x-\tilde{\lambda} \, t}}
        { \mathlarger{t}^{{\scriptscriptstyle \frac{1}{p+1}}} }}_{-\infty}
             \Bigl( \left( \, 
                A-B \xi^2 \, 
                \right)\vee 0 \, \Bigr)^{\frac{1}{p-1}}
                \, \mathrm{d}\xi }
\end{align}
which 
is said to be 
``contact wave for $p$-Laplacian type viscosity''. 
%
Now we are ready to state our main result. 

\medskip

\noindent
{\bf Theorem 1.1} (Main Theorem){\bf .}\quad{\it
Let the flux function $f$ satisfy {\rm(1.4)} and 
the far field states $u_{\pm }$ {\rm(1.5)}. 
Assume that the initial data satisfies 
$u_0-u_0 ^{\rm{R}} \in L^2$ and 
$\partial _x u_0 \in L^{p+1}$. 
Then the Cauchy problem {\rm(1.1)} with $p>1$ has a 
unique global 
solution in time $u=u(t,x)$ 
satisfying 
\begin{eqnarray*}
\left\{\begin{array}{ll}
u-u_0 ^{\rm{R}} \in C^0\bigl( \, [\, 0,\infty)\, ;L^2 \bigr)
                \cap L^{\infty}\bigl( \, \mathbb{R}^{+} \, ;L^{2} \bigr),\\[5pt]
\partial _x u \in L^{\infty} \bigl( \, \mathbb{R}^{+} \, ;L^{p+1} \bigr),\\[5pt]
\partial _t u 
\in L^{\infty} \bigl( \, \mathbb{R}^{+} \, ;L^{p+1} \bigr)
,\\[5pt]
\partial_x \left( \, \left| \, \partial_x u \, \right|^{p-1} \partial_xu \, \right)
\in L^{2}\bigl(\, {\mathbb{R}^{+}_{t}} \times {\mathbb{R}}_{x} \bigr),
\end{array} 
\right.\,
\end{eqnarray*}
and the asymptotic behavior 
$$
\lim _{t \to \infty}\sup_{x\in \mathbb{R}}
|\,u(t,x)-U_{multi}(\, t,x \: ;\: u_-,u_+)\,| = 0, 
$$
where $U_{multi}(t,x)=U_{multi}(\, t,x \: ;\: u_-,u_+)$ is defined as follows: 
in the case $a<u_-<b<u_+$,
$$
U_{multi}(t,x):=
U\left(\, \frac{x-\tilde{\lambda} \, t}{t^{\frac{1}{p+1}}}\: ;\: u_- , b \right)
 + u^r\left(\, \frac{x}{t}\: ;\: b , u_+ \right) - b
$$
and, in the case $u_-<a<b<u_+$,
$$
U_{multi}(t,x):=
u^r\left(\, \frac{x}{t}\: ;\: u_- ,a \right) -a
 + U\left(\, \frac{x-\tilde{\lambda} \, t}{t^{\frac{1}{p+1}}}\: ;\: a ,b \right)
 + u^r\left(\, \frac{x}{t}\: ;\: b , u_+ \right) - b.
$$
}

\medskip

The main theorem is proved by using a technical energy method 
with the aid of the maximum principle, 
and the careful estimates of the nonlinear interactions between 
the nonlinear waves, that is, 
the rarefaction waves and 
the contact wave for the $p$-Laplacian type viscosity. 


This paper is organized as follows. 
In Section 2, we shall prepare the basic properties of 
the rarefaction wave and the 
contact wave for $p$-Laplacian type viscosity. 
In Section 3, we reduce the problem to 
an essential case (similarly in \cite{matsumura-yoshida}, \cite{yoshida})
and reformulate the problem in terms of the deviation from 
the asymptotic state, that is, 
the superposition of the nonlinear waves. 
Following the arguments in \cite{matsu-nishi2}, 
we show 
the global existence of the solution to the reformulated problem 
and the energy estimates which are depending on the time. 
In order to show the asymptotics, 
in Section 4 and Section 5, 
we establish the uniform energy estimates in time 
by using a very technical energy method 
and careful estimates of the interactions 
between the nonlinear waves. 
Finally, in Section 6, 
we prove the asymptotic behavior 
by utilizing the uniform energy estimates 
in Section 4 and Section 5. 
\smallskip

{\bf Some Notation.}\quad 
We denote by $C$ generic positive constants unless 
they need to be distinguished. 
In particular, use $C(\alpha, \beta, \cdots )$ 
or $C_{\alpha, \beta, \cdots }$ 
when we emphasize the dependency on $\alpha, \beta, \cdots $, 
and $\mathbb{R}^{+}$ as 
$
\mathbb{R}^{+}:=(0,\infty).
$
We also use the Friedrichs mollifier $\rho_\delta \ast $, 
where, 
$\rho_\delta(x):=\frac{1}{\delta}\rho \left( \frac{x}{\delta}\right)$
with 
\begin{align*}
\begin{aligned}
&\rho \in C^{\infty}_0(\mathbb{R}),\quad 
\rho (x)\geq 0\:  (x \in \mathbb{R}), \\
&\mathrm{supp} \{\rho \} \subset 
\left\{x \in \mathbb{R}\: \left|\:  |\, x \, |\le 1 \right. \right\},\quad  
\int ^{\infty}_{-\infty} \rho (x)\, \mathrm{d}x=1, 
\end{aligned}
\end{align*}
and $\rho_\delta \ast f$ denote the convolution. 
For function spaces, 
$L^p = L^p(\mathbb{R})$ and $H^k = H^k(\mathbb{R})$ 
denote the usual Lebesgue space and 
$k$-th order Sobolev space on the whole space $\mathbb{R}$ 
with norms $||\cdot||_{L^p}$ and $||\cdot||_{H^k}$, respectively. 
We also define 
the bounded $C^{m}$-class $\mathscr{B}^{m}$ as follows 
$$
f\in \mathscr{B}^{m}(\Omega)
\, \Longleftrightarrow  \, 
f\in C^m(\Omega), 
\; 
\sup _{\Omega}\, \sum _{k=0}^{m} \, \bigl| \, D^kf\, \bigr|<\infty 
$$
for $m< \infty$ and 
$$
f\in \mathscr{B}^{\infty }(\Omega)
\, \Longleftrightarrow  \, 
\forall n\in \mathbb{N},\, f\in C^n(\Omega), 
\; 
\sup _{\Omega}\, \sum _{k=0}^{n} \, \bigl| \, D^kf\, \bigr|<\infty 
$$
where $\Omega \subset \mathbb{R}^d$ and 
$D^k$ denote the all of $k$-th order derivatives. 

\section{Preliminaries} 
In this section, 
we shall arrange the several lemmas concerning with 
the basic properties of 
the rarefaction wave and the viscous contact wave 
for accomplishing the proof of the main theorem. 
Since the rarefaction wave $u^r$ is not smooth enough, 
we need some smooth approximated one 
as in the previous works in 
\cite{hashimoto-matsumura}, \cite{liu-matsumura-nishihara}, \cite{matsu-nishi1}, \cite{matsumura-yoshida}. 
We start with the well-known arguments on $u^r$ 
and the method of constructing its smooth approximation. 
We first consider the rarefaction wave solution $w^r$ 
to the Riemann problem 
for the non-viscous Burgers equation 
\begin{equation}
\label{riemann-burgers}
  \left\{\begin{array}{l}
  \partial _t w + 
  \displaystyle{ \partial _x \left( \, \frac{1}{2} \, w^2 \right) } = 0 
  \, \, \; \; \qquad \quad \qquad ( t > 0,\,x\in \mathbb{R}),\\[7pt]
  w(0,x) = w_0 ^{\rm{R}} ( \, x\: ;\: w_- ,w_+):= \left\{\begin{array}{ll}
                                                  w_+ & \, \: \; \quad (x>0),\\[5pt]
                                                  w_- & \, \: \; \quad (x<0),
                                                  \end{array}
                                                  \right.
  \end{array}
  \right.
\end{equation}
where $w_\pm \in \mathbb{R} \: (w_-<w_+)$ are 
the prescribed far field states. 
The unique global weak solution 
$w=w^r\left( \, \frac{x}{t}\: ;\: w_-,w_+\right)$ 
of (\ref{riemann-burgers}) is explicitly given by 
\begin{equation}
\label{rarefaction-burgers}
w^r \left( \, \frac{x}{t}\: ;\: w_-,w_+\right) := 
  \left\{\begin{array}{ll}
  w_{-} & \bigl(\, x \leq w_{-}t \, \bigr),\\[5pt]
  \displaystyle{ \frac{x}{t} } & \bigl(\, w_{-} t \leq x \leq w_{+} t \, \bigr),\\[5pt]
  w_+ & \bigl(\, x\geq w_{+} t \, \bigr).
  \end{array}\right.
\end{equation} 
Next, under the condition 
$f''(u)>0\ (u\in \mathbb{R})$ and $u_-<u_+$, 
the rarefaction wave solution 
$u=u^r\left( \, \frac{x}{t}\: ;\: u_-,u_+\right)$ 
of the Riemann problem (1.2) 
for hyperbolic conservation law 
is exactly given by 
\begin{equation}
u^r\left( \, \frac{x}{t} \: ; \:  u_-,u_+\right) 
= (\lambda)^{-1}\biggl( w^r\left( \, \frac{x}{t} \: ; \:  \lambda_-,\lambda_+\right)\biggr)
\end{equation}
which is nothing but (1.6), 
where $\lambda_\pm := \lambda(u_\pm) = f'(u_\pm)$. 
We define a smooth approximation of $w^r(\, \frac{x}{t}\: ;\: w_-,w_+)$ 
by the unique classical solution 
$$
w=w(\, t,x\: ;\: w_-,w_+)\in \mathscr{B}^{\infty }( \, [\, 0,\infty )\times \mathbb{R})
$$
to the Cauchy problem for the following 
non-viscous Burgers equation
\begin{eqnarray}
\label{smoothappm}
\left\{\begin{array}{l}
 \partial _t w 
 + \displaystyle{ \partial _x \left( \, \frac{1}{2} \, w^2 \right) } =0 
 \, \, \; \; \quad \qquad \qquad \qquad \qquad \qquad  
 (\ t>0,\,x\in \mathbb{R}),\\[7pt]
 w(0,x) 
 = w_0(x) 
 := \displaystyle{ \frac{w_-+w_+}{2} + \frac{w_+-w_-}{2}\tanh x } 
 \qquad \quad \; \:  (x\in \mathbb{R}),
\end{array}
\right.
\end{eqnarray}   
By using the method of characteristics, 
we get the following formula
\begin{eqnarray}
 \left\{\begin{array} {l}
 w(t,x)=w_0\bigl(x_0(t,x)\bigr)=
 \displaystyle{ \frac{\lambda_-+\lambda_+}{2} } 
+ \displaystyle{ \frac{\lambda_+-\lambda_-}{2}\tanh \bigl( x_0(t,x) \bigr)} ,\\[7pt]
 x=x_0(t,x)+w_0\bigl(x_0(t,x)\bigr)\,t.
 \end{array}
  \right.\,
\end{eqnarray}
We also note the assumption of the flux function $f$ to be 
$\lambda'(u)\left( =\frac{\mathrm{d}^2f}{\mathrm{d}u^2}(u)\right)>0$. 

Now we summarize the results for the smooth approximation $w(\, t,x\: ;\: w_-,w_+)$ 
in the next lemma. 
Since the proof is given by the direct calculation as in \cite{matsu-nishi1}, 
we omit it. 

\medskip

\noindent
{\bf Lemma 2.1.}\quad{\it
Assume that the far field states satisfy $w_-<w_+$. 
Then the classical solution $w(t,x)=w(\, t,x\: ;\: w_-,w_+)$
given by {\rm(2.4)} 
satisfies the following properties: 

\noindent
{\rm (1)}\ \ $w_- < w(t,x) < w_+$ and\ \ $\partial_xw(t,x) > 0$  
\quad  $(t>0, x\in \mathbb{R})$.

\smallskip

\noindent
{\rm (2)}\ For any $1\leq q \leq \infty$, there exists a positive 
constant $C_q$ such that
             \begin{eqnarray*}
                 \begin{array}{l}
                    \parallel \partial_x w(t)\parallel_{L^q} \leq 
                    C_q (1+t)^{-1+\frac{1}{q}} 
                    \; \quad \bigl(t\ge 0 \bigr),\\[5pt]
                    \parallel \partial_x^2 w(t) \parallel_{L^q} \leq 
                    C_q (1+t)^{-1} 
                    \, \, \quad \quad \bigl(t\ge 0 \bigr).
                    \end{array}       
              \end{eqnarray*}
              
\smallskip

\noindent
{\rm (3)}\; $\displaystyle{\lim_{t\to \infty} 
\sup_{x\in \mathbb{R}}
\left| \,w(t,x)- w^r \left( \frac{x}{t} \right) \, \right| = 0}.$
}

\bigskip

\noindent
We define the approximation for 
the rarefaction wave $u^r\left( \, \frac{x}{t}\: ;\: u_-,u_+\right)$ by 
\begin{equation}
U^r(\, t,x\: ; \: u_-,u_+) := (\lambda)^{-1} \bigl( w(\, t,x\: ;\: \lambda_-,\lambda_+)\bigr).
\end{equation}

Then we have the next lemma 
as in the previous works 
(cf. \cite{hashimoto-matsumura}, \cite{liu-matsumura-nishihara}, \cite{matsu-nishi1}, \cite{matsumura-yoshida}). 

\medskip

\noindent
{\bf Lemma 2.2.}\quad{\it
Assume that the far field states satisfy $u_-<u_+$, 
and the flux fanction $f\in C^3(\mathbb{R})$, $f''(u)>0 \: (u\in [\,u_-,u_+\,])$. 
Then we have the following properties:

\noindent
{\rm (1)}\ $U^r(t,x)$ defined by {\rm (2.6)} is 
the unique $C^2$-global solution in space-time 
of the Cauchy problem
$$
\left\{
\begin{array}{l} 
\partial _t U^r +\partial _x \bigl( f(U^r ) \bigr) = 0 
\, \, \, \, \; \; \; \quad \quad \qquad \qquad \qquad \qquad  
(t>0, x\in \mathbb{R}),\\[7pt]
U^r(0,x) 
= \displaystyle{ (\lambda)^{-1} \left( \, \frac{\lambda_- + \lambda_+}{2} 
+ \frac{\lambda_+ - \lambda_-}{2} \tanh x \, \right) } 
\; \: \: \quad \quad( x\in \mathbb{R}),\\[7pt]
\displaystyle{\lim_{x\to \pm \infty}} U^r(t,x) =u_{\pm} 
\, \, \: \: \; \; \quad \quad \qquad \qquad \qquad \qquad \qquad \qquad 
\bigl(t\ge 0 \bigr).
\end{array}
\right.\,     
$$
{\rm (2)}\ \ $u_- < U^r(t,x) < u_+$ and\ \ $\partial_xU^r(t,x) > 0$  
\quad  $(t>0, x\in \mathbb{R})$.

\smallskip

\noindent
{\rm (3)}\ For any $1\leq q \leq \infty$, there exists a positive 
constant $C_q$ such that
             \begin{eqnarray*}
                 \begin{array}{l}
                    \parallel \partial_x 
                    U^r(t) 
                    \parallel_{L^q} \leq 
                    C_q(1+t)^{-1+\frac{1}{q}} 
                    \quad \bigl(t\ge 0 \bigr),\\[5pt]
                    \parallel \partial_x^2 U^r(t) \parallel_{L^q} \leq 
                    C_q(1+t)^{-1}
                    \, \, \quad \quad \bigl(t\ge 0 \bigr).
                    \end{array}       
              \end{eqnarray*}
              
\smallskip

\noindent
{\rm (4)}\; $\displaystyle{\lim_{t\to \infty} 
\sup_{x\in \mathbb{R}}
\left| \,U^r(t,x)- u^r \left( \frac{x}{t} \right) \, \right| = 0}.$

\smallskip

\noindent
{\rm (5)}\ For any $\epsilon \in (0,1)$, there exists a positive 
constant $C_\epsilon$ such that
$$
\left|\, 
U^r(t,x)-u_+ \, \right|
\leq C_\epsilon (1+t)^{-1+\epsilon}
     \mathrm{e}^{-\epsilon \, | \, x-\lambda_+t \, |}
\quad \bigl(t\ge 0, x \ge \lambda_+t \bigr).
$$

\smallskip

\noindent
{\rm (6)}\ For any $\epsilon \in (0,1)$, there exists a positive 
constant $C_\epsilon$ such that
$$
\left|\, 
U^r(t,x)-u_- \, \right|
\leq C_\epsilon(1+t)^{-1+\epsilon}
     \mathrm{e}^{-\epsilon \, | \, x-\lambda_-t \, |}
\quad \bigl(t\ge 0, x \le \lambda_-t\bigr).
$$

\smallskip

\noindent
{\rm (7)}\ For any $\epsilon \in (0,1)$, there exists a positive 
constant $C_\epsilon$ such that
$$
\left| \,U^r(t,x) - u^r\left( \frac{x}{t}\right) \, \right| 
\leq C_\epsilon(1+t)^{-1+\epsilon} 
\qquad \bigl(t \ge 1,  \lambda_-t \le x \le \lambda_+t\bigr).
$$

\smallskip

\noindent
{\rm (8)}\ For any $(\epsilon ,q)\in (0,1)\times [\, 1,\infty \, ]$, 
there exists a positive 
constant $C_{\epsilon,q}$ such that
$$
\left|\left|
 \,U^r(t,\cdot \: ) - u^r\left( \frac{\cdot }{t}\right) \, 
\right|\right|_{L^q}  
\leq C_{\epsilon,q}(1+t)^{-1+\frac{1}{q}+\epsilon} 
\qquad \bigl(t \ge 0 \bigr).
$$
}

\noindent
Because the proofs of (1) to (4) are given in \cite{matsu-nishi1}, 
(5) to (7) are in \cite{matsumura-yoshida} 
and (8) is in \cite{yoshida}, 
we omit the proofs here. 
\bigskip

We also prepare the next lemma for the properties of 
the contact wave for $p$-Laplacian type viscosity 
$U\Bigl(\, \frac{x}{t^{\frac{1}{p+1}}}\,;\,u_- ,u_+ \Bigr)$ defined by (1.11). 
In the following, we abbreviate 
``contact wave for $p$-Laplacian type viscosity'' 
to ``viscous contact wave''.  
Substituting (1.12) into (1.13), 
we rewrite the viscous contact wave as 
\begin{align}
\begin{aligned}
U(t,x)&=U\left(\, \frac{x}{t^{\frac{1}{p+1}}}\: ;\: u_- ,u_+ \right) \\
      &=u_{+}-\Lint^{\infty}_{x}
             \frac{1}{t^{\frac{1}{p+1}}}\, 
             \bigggl( \left( 
             A-B \left(\frac{y}{t^{\frac{1}{p+1}}} \right)^2 \, 
             \right)\vee 0 \bigggr)^{\frac{1}{p-1}}
             \, \mathrm{d}y, 
\end{aligned}
\end{align}
where 
\begin{eqnarray*}
\left\{\begin{array}{ll}
\displaystyle{
A=A_{p,\mu,u_{\pm}}
:= \left( \, 
   \frac{(p-1)\, \left( \, u_{+} - u_{-} \, \right)}
   {8\, \mu\, p(p+1) 
   \biggl( \, \displaystyle{\int^{\frac{\pi}{2}}_{0} 
   \bigl( \sin \theta \bigr)^{\frac{p+1}{p-1}}\, \mathrm{d}\theta }
   \biggr)^2 }
   \, \right)^{\frac{p-1}{p+1}}
}, \\[25pt]
\displaystyle{ 
B=B_{p,\mu}:= \frac{p-1}{2\, \mu \, p(p+1)} 
}, \\[10pt] 
\displaystyle{ 
2A^{\frac{p+1}{2(p-1)}}B^{-\frac{1}{2}}
      \int^{\frac{\pi}{2}}_{0} 
      \bigl( \sin \theta \bigr)^{\frac{p+1}{p-1}}\, \mathrm{d}\theta 
      = u_{+} - u_{-}}
.
\end{array}
\right.\,
\end{eqnarray*}

\noindent
Then, we have the next lemma. 
Because the proofs are very elementary, 
we omit the proofs. 

\medskip

\noindent
{\bf Lemma 2.3.}\quad{\it
For any 
$p>1$ and $u_\pm \in \mathbb{R}$, we have the following: 

\noindent
{\rm (i)}\ $U$ defined by {\rm (1.11)} satisfies 
$$
U \in 
\mathscr{B}^1\bigl( \, (0,\infty \bigr)\times \mathbb{R}\, ) 
\mathlarger{ \mathlarger{ \mathlarger{ 
\setminus } } } \, 
C^2
\left( \, 
\biggl\{ \, 
(t,x) \in \mathbb{R}^{+} \times \mathbb{R} \, \biggr|\biggl. \, 
x = \pm \sqrt{\frac{A}{B}}\, t^{\frac{1}{p+1}} \, 
\biggr\} \, \right), 
$$  
and is a 
self-similar type strong solution of the Cauchy problem 
   $$
           \left\{
              \begin{array}{l} 
              \partial _t U 
              - \mu \, \partial_x \left( \, 
              \left| \, \partial_x U \, \right|^{p-1} \partial_x U \, \right) 
              = 0 
              \quad \qquad \qquad  (t>0, x\in \mathbb{R}),\\[13pt]
              U(0,x) = u_0 ^{\rm{R}} (\, x\: ;\: u_- ,u_+)
               = 
               \left\{\begin{array} {ll}
               u_-  & \, \; \; \quad \qquad (x < 0),\\[5pt]
               u_+  & \, \; \; \quad \qquad (x > 0),
               \end{array}\right.\\[13pt]
              \displaystyle{\lim_{x\to \pm \infty}} U(t,x) =u_{\pm} 
              \, \: \: \: \; \quad \quad \qquad \qquad \qquad \qquad \bigl(\, t\ge 0\bigr).
              \end{array}
            \right.\,   
    $$          
{\rm (ii)}\ For $t>0$ and $x\in \mathbb{R}$, 
$$
\left\{
   \begin{array}{l} 
   U(t,x)=u_{-},
   \, \: \; \quad \qquad \qquad \qquad \qquad 
   \left( \, x \leq -\sqrt{\frac{A}{B}}\, t^{\frac{1}{p+1}} \,  \right),
   \\[5pt]
   u_- < U(t,x) < u_+,\; \partial_xU(t,x) > 0 \quad 
   \left( \, 
   -\sqrt{\frac{A}{B}}\, t^{\frac{1}{p+1}} < x < \sqrt{\frac{A}{B}}\, t^{\frac{1}{p+1}}
   \, \right),
   \\[5pt]
   U(t,x)=u_{+},
   \, \: \; \quad \qquad \qquad \qquad \qquad 
   \left( \,  x \geq \sqrt{\frac{A}{B}}\, t^{\frac{1}{p+1}} \,  \right). 
   \end{array}
\right.\, 
$$

\smallskip

\noindent
{\rm (iii)}\ It holds that for any $1 \leq  q < \infty $, 
$$
\|\, \partial _x 
U (t)\, \| _{L^q} 
= C_1(\, A,B\: ;\: p ,q\, ) \, t^{-\frac{q-1}{(p+1)q}}\qquad (t > 0) 
$$
where 
$$
C_1(\, A,B\: ;\: p ,q\, ):=
\left( 
2A^{\frac{p+2q-1}{2(p-1)}}B^{-\frac{1}{2}}
\int^{\frac{\pi}{2}}_{0} 
      \bigl( \sin \theta \bigr)^{\frac{q}{p-1}}
\, \mathrm{d}\theta 
\right)^{\frac{1}{q}}.
$$
If $q = \infty$, we have 
$$
\|\, \partial _x U (t) \, \| _{L^\infty} 
= \left(2A \right)^{\frac{1}{p-1}} \, t^{- \frac{1}{p+1}}\qquad (t > 0).
$$

\smallskip

\noindent
{\rm (iv)}\ It holds that for any $1 \leq  q < \frac{p-1}{p-2}$ with $p>2$, 
or any $1 \leq  q < \infty$ with $1 < p \leq 2$, 
$$
\|\, \partial _x^2 U (t) \, \| _{L^q} 
= C_2(\, A,B\: ;\: p ,q\, ) \, t^{-\frac{2q-1}{(p+1)q}}\qquad (t > 0) 
$$
where 
\begin{align*}
\begin{aligned}
& C_2(\, A,B\: ;\: p ,q\, ) \\
& \quad \, 
:=
\bigggl( 
2\left(
\frac{2A^{-\frac{ p-2 }{p-1}} B}{p-1}
\right)^{q}
\left( \frac{B}{A} \right)^{- \frac{q+1}{2}}
\int^{\frac{\pi}{2}}_{0} 
      \bigl( \sin \theta \bigr)^{- \frac{2 ( p-2 )q}{p-1} + 1}
      \bigl( \cos \theta \bigr)^{q}
\, \mathrm{d}\theta
\bigggr)^{\frac{1}{q}}. 
\end{aligned}
\end{align*}
If $1 < p \leq 2$, for $q = \infty$, we have 
$$
\|\, \partial _x^2 U (t) \, \| _{L^\infty} 
= {\frac{2A^{\frac{\left| p-2 \right|}{p-1}} B}{p-1}} 
     \left( \frac{B}{A} \right)^{-\frac{1}{2}}\, 
     t^{- \frac{2}{p+1}}\qquad (t > 0).
$$

\smallskip

\noindent
{\rm (v)}\ It holds that 
$$
\left|\left| \, 
\partial_x \left( \, 
\left| \, \partial_x U \, \right|^{p-1} \partial_x U \, 
\right) (t) \, 
\right| \right| _{L^2} 
= C_3(\, A,B\: ;\: p \, ) \, t^{-\frac{2p+1}{2(p+1)}}\qquad (t > 0) 
$$
where 
$$
C_3(\, A,B\: ;\: p \, ):=
\left( 
2\left(
\frac{2B^{p}}{p-1}
\right)^{2}
\left( \frac{B}{A} \right)^{-\frac{3p-7}{2(p-1)}}
\int^{\frac{\pi}{2}}_{0} 
      \bigl( \sin \theta \bigr)^{\frac{p+3}{p-1}}
      \bigl( \cos \theta \bigr)^{2}
\, \mathrm{d}\theta
\right)^{\frac{1}{2}}.
$$
\smallskip

\noindent
{\rm (vi)}\; $\displaystyle{\lim_{t\to \infty} \sup_{x\in \mathbb{R}}\, 
\bigl|\,U(1+t,x)- U(t,x)\, \bigr| = 0}.$
%

}
%
%
%

\section{Reformulation of the problem}
In this section, 
we reduce our Cauchy problem (1.1) to 
a simpler case 
and reformulate the problem 
in terms of the deviation from the asymptotic state 
(the same as in \cite{matsumura-yoshida}, \cite{yoshida}). 
At first, 
without loss of generality, 
we shall consider the case 
where $a<0$, $b=0$ and the flux function $f(u)$ satisfies 
\begin{equation}
\left\{
\begin{array}{ll}
  f''(u) >0 & \; (u \in (-\infty ,a]\cup [0,+\infty )),\\[5pt]
  f(u) =0 & \; (u \in (a,0)),
\end{array}\right.
\end{equation}
under changing the variables and constant as
$x-\tilde{\lambda} \,t \mapsto x$,  $u-b \mapsto u$,
$f(u+b)-f'(b)\,u-f(a) \mapsto f(u)$ and $a-b \mapsto a$
in this order. 
For the far field states 
$u_\pm \in \mathbb{R}$, 
we only deal with the typical case 
$a<u_-<0<u_+$ for simplicity, since the case 
$u_-<a<0<u_+$ can be treated 
technically 
in the same way of the proof as 
$a<u_-<0<u_+$. 
Indeed, in the case $u_-<a<0<u_+$,
as we shall see in Section 4  and Section 5, 
there appears the extra nonlinear interaction terms
between two rarefaction waves 
$u^r(\, \frac{x}{t}\: ;\: u_-,a)$ and $u^r(\, \frac{x}{t}\: ;\: 0,u_+)$
with $\lambda(a)=\lambda(0)=0$ in the remainder term of 
the viscous conservation law for the asymptotics 
$U_{multi}$\,(see the right-hand side of (3.5)). 
These terms can be handled in much easier way 
by Lemma 2.2 than that for other essential nonlinear interaction terms 
between the rarefaction and the viscous contact waves. 
Furthermore, we should point out that 
the problem under the assumptions for the flux function (3.1) 
and the far field states $a<u_-<0<u_+$
is essentially the same as that for $a=-\infty$, because 
obtaining the a priori and the uniform energy estimates for the former one can be
given in almost the same way as the latter one. 
Therefore, it is quite natural for us 
to treat only 
a simple case 
\begin{equation}
\left\{
\begin{array}{ll}
  f''(u) >0 & \; (u \in [\, 0,\infty )),\\[5pt]
  f(u) =0 & \; (u \in (-\infty,0)),
\end{array}\right.
\end{equation}
and assume $u_-<0<u_+$.
The corresponding main theorem is the following. 

\medskip

\noindent
{\bf Theorem 3.1.}\quad{\it
Let the flux function 
$f\in C^1(\mathbb{R})\cap C^3 (\, [\, 0,\infty))$  
satisfy {\rm(3.2)} and 
the far field states $u_{-}<0<u_{+}$. 
Assume that the initial data satisfies 
$u_0-u_0 ^{\rm{R}} \in L^2$ and 
$\partial _x u_0 \in L^{p+1}$. 
Then the Cauchy problem {\rm(1.1)} with $p>1$ has a 
unique global solution in time $u$ 
satisfying 
\begin{eqnarray*}
\left\{\begin{array}{ll}
u-u_0 ^{\rm{R}} \in C^0\bigl( \, [\, 0,\infty)\, ;L^2 \bigr)
                    \cap L^{\infty}\bigl( \, \mathbb{R}^{+} \, ;L^{2} \bigr),\\[5pt]
\partial _x u \in L^{\infty} \bigl( \, \mathbb{R}^{+} \, ;L^{p+1} \bigr)
                  \cap L^{p+2}\bigl(\, {\mathbb{R}^{+}_{t}} 
                       \times 
                       \left\{x \in \mathbb{R}\, |\, u>0 \right\} \bigr),\\[5pt]
\partial _t u 
\in  L^{\infty} \bigl( \, \mathbb{R}^{+} \, ;L^{p+1} \bigr)
,\\[5pt]
\partial_x \left( \, \left| \, \partial_xu \, \right|^{p-1} \partial_xu \, \right)
\in L^{2}\bigl(\, {\mathbb{R}^{+}_{t}} \times {\mathbb{R}}_{x} \bigr),
\end{array} 
\right.\,
\end{eqnarray*}
and the asymptotic behavior 
$$
\lim _{t \to \infty}\sup_{x\in \mathbb{R}}
|\,u(t,x)-U_{multi}(\, t,x \: ;\: u_-,u_+)\,| = 0, 
$$
where $U_{multi}(t,x)=U_{multi}(\, t,x \: ;\: u_-,u_+)$ is defined by 
$$
U_{multi}(t,x):=
 U\left(\, \frac{x}{t^{\frac{1}{p+1}}}\: ;\: u_- ,0 \right)
 + u^r\left(\, \frac{x}{t}\: ;\: 0 , u_+ \right).
$$
}

\medskip
Here, we first should note by Lemma 2.2 and Lemma 2.3, 
the asymptotic state $U_{multi}(\, t,x\: ; \: u_-, u_+)$ 
can be replaced by a following approximated one 
\begin{equation}
\tilde{U}(t,x) := U(1+t,x) + U^r(t,x)
\end{equation}
where
$$
U(1+t,x)=U\left( \frac{x}{(1+t)^{\frac{1}{p+1}}}\: ;\: u_- ,0 \right),
\quad U^r(t,x) =U^r(\, t,x\: ; \: 0, u_+).
$$
This is because, 
from Lemma 2.2 (especially (4)) and Lemma 2.3 (especially (vi)), 
\begin{align*}
\begin{aligned}
&\sup_{x\in \mathbb{R}}\, \left| \, \tilde{U}(t,x) - U_{multi}(t,x) \, \right| 
\le
\sup_{x\in \mathbb{R}} \, \bigl| \, U(1+t,x) - U(t,x) \, \bigr|\\ 
&\qquad \qquad + \sup_{x\in \mathbb{R}}\, 
\left| \, U^r(t,x) - u^r\left( \frac{x}{t}\right) \, \right|
\to 0\quad (t \to \infty).
\end{aligned}
\end{align*}
In the following, we write $U(1+t,x)$ again $U(t,x)$ for simplicity. 
Then it is noted that $\tilde{U}$ approximately satisfies 
the equation of (1.1) as
\begin{equation}
 \partial _t\tilde{U} +\partial_x \bigl(f(\tilde{U})\bigr)
 - \mu \, 
 \partial_x 
\left( \, 
\bigl| \, \partial_x \tilde{U} \, \bigr|^{p-1} \partial_x \tilde{U} \, 
\right)
 = -F_{p}(U,U^r),
\end{equation}
where the remainder term $F_{p}(U,U^r)$ is explicitly given by
\begin{align}
\begin{aligned}
F_{p}(U,U^r)
:= 
&\widetilde{F_{p} }(U,U^r)\\
&+ \mu \, \partial_x 
\left( \, 
\bigl| \, \partial_x U + \partial_x U^{r} \, \bigr|^{p-1} 
\bigl( \, \partial_x U + \partial_x U^{r} \, \bigr) - 
\left| \, \partial_x U \, \right|^{p-1} \partial_x U  \, 
\right)\\
:= 
&-\bigl( \, f'(U+U^r)
- f'(U^r) \, \bigr) \, \partial_x U^r 
- f'(U+U^r)\, \partial_x U \\
&+ \mu \, \partial_x 
\left( \, 
\bigl| \, \partial_x U + \partial_x U^{r} \, \bigr|^{p-1} 
\bigl( \, \partial_x U + \partial_x U^{r} \, \bigr) - 
\left| \, \partial_x U \, \right|^{p-1} \partial_x U  \, 
\right)
\end{aligned}
\end{align}
which consists of the interaction terms of the viscous contact wave $U$ and 
the approximation of the rarefaction wave $U^r$, 
and the approximation error of $U^r$ as solution to 
the conservation law for the $p$-Laplacian type viscosity. 
Here we should 
note that 
$U$ is monotonically nondecreasing 
and $U^r$ is monotonically increasing, 
that is, 
$\partial_x\tilde{U}(t,x) > 0\ \bigl(t\ge 0, x\in \mathbb{R} \bigr)$
which is frequently used hereinafter. 
Now 
putting 
\begin{equation}
u(t,x) = \tilde{U}(t,x) + \phi(t,x)
\end{equation}
and using (3.5),
we can reformulate the problem (1.1) in terms of 
the deviation $\phi $ from 
$\tilde{U}$ as 
\begin{eqnarray}
 \left\{\begin{array}{ll}
  \partial _t\phi + \partial_x \left( f(\tilde{U}+\phi) - f(\tilde{U}) \right) \\[5pt]
  \quad - \mu \, \partial_x 
    \left( \, 
    \bigl| \, \partial_x \tilde{U} + \partial_x \phi \, \bigr|^{p-1} 
    \bigl( \, \partial_x \tilde{U} + \partial_x \phi \, \bigr) - 
    \bigl| \, \partial_x \tilde{U} \, \bigr|^{p-1} \partial_x \tilde{U}  \, 
    \right)\\[2pt]
    \quad \qquad \qquad \qquad \qquad \qquad \qquad = F_{p}(U,U^r)
  \, \, \: \: \quad  (t>0, x\in \mathbb{R}), \\[5pt]
  \phi(0,x) = \phi_0(x) 
  := u_0(x)-\tilde{U}(0,x) 
  \qquad \qquad \quad \; \: \: \; \; \: \; \; \,\, \, \, (x\in \mathbb{R}). 
 \end{array}
 \right.\,
\end{eqnarray}
Then we look for 
the unique global solution in time 
$\phi $ which has the asymptotic behavior 
$$
\displaystyle{
\sup_{x \in \mathbb{R}}\left|\,  \phi(t,x) \, \right|
\xrightarrow[\, \, t\to  \infty \, ]
\ \ 0}. 
$$
Here we note the fact $\phi_0 \in L^2$ 
and $\partial_x \phi_0 \in L^{p+1}$ by the assumptions on $u_0$ 
and the fact 
$$
\partial_x \tilde{U}(0,\cdot \, )
=\partial_x U(0,\cdot \, )+\partial_x U^r(0,\cdot \, )\in L^{p+1}.
$$
In the following, we always assume that 
the flux function $f\in C^1(\mathbb{R})\cap C^3 (\, [\, 0,\infty))$ satisfies (3.2), 
and  
the far field states satisfy $u_-<0<u_+$.
Then the corresponding our main theorem for $\phi$ we should prove is 
as follows. 

\medskip

\noindent
{\bf Theorem 3.2.}\quad{\it
Suppose 
$\phi_0 \in L^2$ and 
$\partial _x \phi_0 \in L^{p+1}$. 
Then there exists 
the unique global solution in time $\phi=\phi(t,x)$ 
of the 
Cauchy problem {\rm (3.7)} 
satisfying 
\begin{eqnarray*}
\left\{\begin{array}{ll}
\phi \in C^0\bigl( \, [\, 0,\infty)\, ;L^2 \bigr)
            \cap L^{\infty}\bigl( \, \mathbb{R}^{+} \, ;L^{2} \bigr),\\[5pt]
\partial _x \phi \in L^{\infty} \bigl( \, \mathbb{R}^{+} \, ;L^{p+1} \bigr)
                     \cap L^{p+1}
                          \bigl(\, {\mathbb{R}^{+}_{t}} 
                          \times {\mathbb{R}}_{x} \bigr),\\[5pt]
\partial _x \bigl( \, \tilde{U} + \phi \, \bigr)
              \in L^{\infty} \bigl( \, \mathbb{R}^{+} \, ;L^{p+1} \bigr)
                  \cap L^{p+2}\bigl(\, {\mathbb{R}^{+}_{t}} 
                       \times 
                       \left\{x \in \mathbb{R}\, |\, u>0 \right\} \bigr),\\[5pt]
\partial _t \bigl( \, \tilde{U} + \phi \, \bigr) 
\in  L^{\infty} \bigl( \, \mathbb{R}^{+} \, ;L^{p+1} \bigr)
,\\[5pt]
\partial_x \left( \, \bigl| \, 
\partial_x\bigl( \, \tilde{U} + \phi \, \bigr)\, 
\bigr|^{p-1} \partial_x\bigl( \, \tilde{U} + \phi \, \bigr) \, \right)
\in L^{2}\bigl(\, {\mathbb{R}^{+}_{t}} \times {\mathbb{R}}_{x} \bigr),
\end{array} 
\right.\,
\end{eqnarray*}
and the asymptotic behavior 
$$
\lim _{t \to \infty}\sup_{x\in \mathbb{R}}
|\,\phi (t,x)\,| = 0. 
$$
}

In order to accomplish the proof of Theorem 3.2, 
we first note that for any $T>0$, the global existence on $[\, 0,T\,]$ 
and uniqueness can be proved by the similar arguments as in \cite{matsu-nishi2}. 
Indeed, we rewrite our Cauchy problem {\rm (3.7)} again as 
\begin{eqnarray}
 \left\{\begin{array}{ll}
    \partial _t \phi + \partial_x \bigl( \, f(U+U^{r}+\phi) - f(U^{r}) \, \bigr)  
    + \mu \, \partial_x 
    \left( \, 
    \bigl( \, \partial_x U \, \bigr)^{p} \, 
    \right) \\[5pt]
    \; \; \; \; \; \: \: 
      = \mu \, \partial_x 
      \left( \, 
      \bigl| \, \partial_x (U+U^{r}+\phi) \, \bigr|^{p-1} 
      \partial_x (U+U^{r}+\phi) \,  
      \right) 
  \, \, \: \: \quad  (t>0, x\in \mathbb{R}), \\[10pt]
  \phi(0,x) = \phi_0(x) 
  := u_0(x)-\tilde{U}(0,x)
  \qquad \qquad \; \: \qquad \quad \qquad 
  \; \: \: \; \; \: \; \; \,\, \, \, 
  (x\in \mathbb{R}), 
 \end{array}
 \right.\,
\end{eqnarray} 
and for any $\epsilon \in (0,1 \,]$, 
we consider the $\epsilon$-regularized problem as 
\begin{eqnarray}
 \left\{\begin{array}{ll}
  \partial _t \phi_{\epsilon} + 
  \displaystyle{ 
  \partial_x 
  \Bigl( \, 
  f^{\epsilon} (U^{\epsilon}+U^{r, \epsilon}+\phi_{\epsilon}) 
  - f^{\epsilon}(U^{r, \epsilon}) \, 
  \Bigr) } 
    + \mu \, \partial_x 
    \left( \, 
    \bigl( \, \partial_x U^{\epsilon} \, \bigr)^{p} \, 
    \right)\\[5pt]
    \quad \, \, 
    = \mu \, \partial_x 
      \left( 
      \left( 
      \bigl( \, \partial_x (U^{\epsilon}+U^{r, \epsilon}+\phi_{\epsilon}) \, 
      \bigr)^2 + \epsilon \, 
      \right)^{\frac{p-1}{2}} \partial_x (U^{\epsilon}+U^{r, \epsilon}+\phi_{\epsilon})   
      \right) \\[10pt]
  \: \; \; \; \; \qquad \qquad \qquad \; \: \: \qquad \qquad \qquad \qquad \qquad \qquad 
  \quad \, 
  (t>0, x\in \mathbb{R}), \\[10pt]
  \phi_{\epsilon}(0,x) = \phi_0^{\epsilon}(x) 
  \; \; \; \; \quad \quad \quad \; \: \: \quad \qquad \qquad 
  \qquad \qquad \qquad \quad 
  (x\in \mathbb{R}). 
 \end{array}
 \right.\,
\end{eqnarray}
where, 
\begin{eqnarray*}
 \left\{\begin{array}{ll}
\phi^{\epsilon} _{0}(x) := ( \rho_{\epsilon} \underset{x}{\ast } \phi_{0} ) (x) 
                        := \displaystyle{\int_{-\infty}^{\infty} 
                           \rho_{\epsilon} (x-y) \, 
                           \phi_{0}(y) \, \mathrm{d}y} 
                           \in H^{\infty}(\mathbb{R}_{x}), \\[10pt]
f^{\epsilon} (u) := ( \rho_{\epsilon} \underset{u}{\ast } f ) (u) 
                 := \displaystyle{\int_{-\infty}^{\infty} 
                    \rho_{\epsilon} (u-v) \, 
                    f(v) \, \mathrm{d}v}
                    \in C^{\infty}(\mathbb{R}_{u}), \\[10pt]
U^{r, \epsilon} (t,x) := \left( \bigl( f^{\epsilon} \bigr)' \right)^{-1} 
                         \Bigl( 
                         w\bigl(\, t,x\: 
                         ;\: f^{\epsilon}(0), \, f^{\epsilon}(w_+) \, \bigr) 
                         \Bigr)
                      \in \mathscr{B}^{\infty}
                          \bigl(\, [ \, 0, \infty)_{t} 
                          \times \mathbb{R}_{x} \bigr), \\[10pt]
U^{\epsilon} (t,x) := ( \rho_{\epsilon} \underset{x}{\ast } U ) (t,x) 
                   := \displaystyle{\int_{-\infty}^{\infty} 
                      \rho_{\epsilon} (x-y) \, 
                      U\left( \frac{y}{(1+t)^{\frac{1}{p+1}}} \right)
                      \, \mathrm{d}y} \\[12pt]
                      \qquad \qquad \qquad \qquad \quad \quad \quad 
                      \in \mathscr{B}^{1} \cap C^{\infty}
                          \bigl(\, (-1, \infty)_{t} \times \mathbb{R}_{x} \bigr).
\end{array}
 \right.\,
\end{eqnarray*}
Here, $w(t,x)$ 
is the classical solution of (2.4) with $w_- = 0$. 
If we define 
$u_{\epsilon} := \phi_{\epsilon} + U^{\epsilon} +U^{r, \epsilon}$, 
we also have the equivalent form 
\begin{eqnarray}
 \left\{\begin{array}{ll}
  \partial _t u_{\epsilon} + 
  \displaystyle{ 
  \partial_x 
  \Bigl( \, 
  f^{\epsilon} (u_{\epsilon}) \, 
  \Bigr) } 
    = \mu \, \partial_x 
      \left( 
      \left( 
      \bigl( \, \partial_x u_{\epsilon} \, 
      \bigr)^2 + \epsilon \, 
      \right)^{\frac{p-1}{2}} \partial_x u_{\epsilon} \, 
      \right) \\[10pt]
  \, \qquad \qquad \qquad 
  \qquad \qquad \qquad \qquad \qquad \qquad 
  \quad \, 
  (t>0, x\in \mathbb{R}), \\[10pt]
  u_{\epsilon}(0,x) = u_0^{\epsilon}(x) 
                    := \phi_0^{\epsilon}(x) 
                       + U^{\epsilon}(0,x)
                       + U^{r, \epsilon}(0,x)
  \quad \; \: \, \, 
  (x\in \mathbb{R}). 
 \end{array}
 \right.\,
\end{eqnarray}
Owing to Lady{\v{z}}enskaja-Solonnikov-Ural'ceva \cite{lad-sol-ura} 
(see also \cite{lions}), 
the regularized problem (3.9) has a unique classical global solution in time 
$\phi_{\epsilon}=\phi_{\epsilon}(t,x)$ 
on $[\, 0,T\, ]\times \mathbb{R} $ for any $T>0$, that is, 
$$
\phi_{\epsilon},\, \partial_x \phi_{\epsilon},\, 
\partial_x^2 \phi_{\epsilon},\, \partial_t \phi_{\epsilon}
\in C^{\infty }\bigl( \, [\, 0,\infty )\times \mathbb{R} \, \bigr)
$$
because the equation in (3.9) is uniformly parabolic 
with variable coefficients. 
Further, 
the maximum principles 
(see \cite{ilin-kalashnikov-oleinik}, \cite{pro-wei}) for (3.10) allows us to get 
the uniform boundedness to $\phi_{\epsilon}(t,x)$ as follows. 

\medskip

\noindent
{\bf Lemma 3.1} (uniform boundedness){\bf .}\quad {\it
It holds that 
\begin{align*}
\begin{aligned}
&\sup_{t\in [\, 0,\infty ), x\in \mathbb{R}}|\,\phi_{\epsilon}(t,x)\,| \\
&\le 
\sup_{x\in \mathbb{R}}|\,u_0^{\epsilon} (x)\,| 
+ \sup_{t\in [\, 0,\infty ), x\in \mathbb{R}}|\,U (t,x)\,| 
+ \sup_{t\in [\, 0,\infty ), x\in \mathbb{R}}|\,U^r (t,x)\,| \, 
 \\ 
&= 
\Vert \,\phi_0 \, \Vert_{L^{\infty}} + 2\, |\, u_{-} \, | + 2\, |\, u_{+} \, | \, 
=: \widetilde{C}. 
\end{aligned}
\end{align*}
}

\medskip

\noindent
Since $\phi^{\epsilon} _{0} \in H^{\infty}$, 
by using the above uniform boundedness and the standard arguments in the Sobolev space 
on the uniformly parabolic equations, 
we can see the classical $C^{\infty}$-solution $\phi_{\epsilon}$ 
also satisfies 
$$
\phi_{\epsilon} 
\in C^{\infty }\bigl( \, [\, 0,\infty )\: ; \: H^{\infty} \bigr). 
$$
Then, for any fixed $\epsilon \in (0,1)$ and $T>0$, 
we can obtain the following a priori estimates 
which are depend upon 
$\epsilon$ and $T$ to the problem (3.9) as follows. 

\medskip

\noindent
{\bf Lemma 3.2} (a priori estimates I){\bf .}\quad {\it
There exists a positive constant $C_{I}$ such that for $0<t<T$, 
\begin{align*}
\begin{aligned}
&\,\Vert \, \phi_{\epsilon}
 (t) \, \Vert_{L^2}^2 
 + \int ^{t }_{0 } \int ^{\infty}_{-\infty} 
 {\phi}^2 \, 
 \bigl( \, \partial _x U^{\epsilon} + \partial _x U^{r, \epsilon} \, \bigr) \, 
 \mathrm{d}x \mathrm{d}\tau  \\
&\qquad \qquad \; \, 
 + \int ^{t }_{0 }\int ^{\infty }_{-\infty } 
   \bigl( \, \partial_x u_{\epsilon} \, \bigr)^2
      \left( 
      \bigl( \, \partial_x u_{\epsilon} \, 
      \bigr)^2 + \epsilon \, 
      \right)^{\frac{p-1}{2}} \, 
  \mathrm{d}x \mathrm{d}\tau 
\leq C_{I},  
\end{aligned}
\end{align*}
where 
\begin{align*}
\begin{aligned}
C_{I}
=C
 \biggl( \, 
 T, \: 
 \widetilde{C}, \: 
 \| \, \phi_0 \, \| _{L^2} 
 \, \biggr).
\end{aligned}
\end{align*}
}

\medskip

\noindent
{\bf Lemma 3.3} (a priori estimates I\hspace{-.1em}I){\bf .}\quad {\it
There exists a positive constant $C_{I\hspace{-.1em}I}$ such that for $0<t<T$, 
\begin{align*}
\begin{aligned}
&\,\Vert \, \partial_x u_{\epsilon}(t) \, \Vert_{L^{p+1}}^{p+1} 
 + {\epsilon}^{\frac{p-1}{2}}
   \,\Vert \, \partial_x u_{\epsilon}(t) \, \Vert_{L^{2}}^{2} \\
&\qquad \qquad \quad \quad \: \, 
 + \int ^{t }_{0 }\int ^{\infty }_{-\infty } 
      \left( 
      \bigl( \, \partial_x u_{\epsilon} \, 
      \bigr)^2 + \epsilon \, 
      \right)^{p-1} 
      \bigl( \, \partial_x^2 u_{\epsilon} \, 
      \bigr)^2 \, 
   \mathrm{d}x \mathrm{d}\tau \\
&\qquad \qquad \quad \quad \: \, 
 + \int ^{t }_{0 }\int ^{\infty }_{-\infty } 
      \left( 
      \bigl( \, \partial_x u_{\epsilon} \, 
      \bigr)^2 + \epsilon \, 
      \right)^{\frac{p-1}{2}} 
      \bigl| \, \partial_x u_{\epsilon} \, 
      \bigr|^3 \, 
   \mathrm{d}x \mathrm{d}\tau 
\leq C_{I\hspace{-.1em}I},  
\end{aligned}
\end{align*}
where 
\begin{align*}
\begin{aligned}
&C_{I\hspace{-.1em}I}
=C
 \biggl( \, 
 T,  \: 
 \widetilde{C}, \: \biggr. \: 
 \| \, \phi_0 \, \| _{L^2}, \: 
\| \, 
\partial_x u_0 
\, \| _{L^{p+1}}, \: 
{\epsilon}^{\frac{p-1}{2}} \, \| \, 
\partial_x u_0 
\, \| _{L^{2}} 
 \, \biggr).
\end{aligned}
\end{align*}
}

\medskip

\noindent
The proofs of Lemma 3.2 and Lemma 3.3 are given in the last part in this section. 
We can also prove the following lemma. 
Because the proof is given in the same way as the above lemma, 
we omit the proof. 
%

\medskip

\noindent
{\bf Lemma 3.4} (a priori estimates I\hspace{-.1em}I\hspace{-.1em}I){\bf .}\quad {\it
There exists a positive constant $C_{I\hspace{-.1em}I\hspace{-.1em}I}$ 
such that for $0<t<T$, 
\begin{align*}
\begin{aligned}
\,\Vert \, \partial_x u_{\epsilon}
(t) \, \Vert_{L^{2}} 
\leq C_{I\hspace{-.1em}I\hspace{-.1em}I}
\end{aligned}
\end{align*}
where 
\begin{align*}
\begin{aligned}
C_{I\hspace{-.1em}I\hspace{-.1em}I}
=C
 \biggl( \, 
 T, \: 
 \widetilde{C}, \: 
 \| \, \phi_0 \, \| _{L^2}, \: 
\| \, 
\partial_x u_0 
\, \| _{L^{p+1}}, \: 
\| \, 
\partial_x u_0 
\, \| _{L^{2}} 
 \, \biggr).
\end{aligned}
\end{align*}
}
%
%
Once the Lemmas 3.2, 3.3 and 3.4 are proved, 
taking the limit $\epsilon \searrow 0$ as in the arguments 
in \cite{matsu-nishi2}, 
we have Theorem 3.3. 

\medskip

\noindent
{\bf Theorem 3.3.}\quad{\it
For any initial data 
$\phi_0 \in L^2$ and 
$\partial _x \phi_0 \in L^{p+1}$, 
there exists the unique global solution in time 
of the 
Cauchy problem {\rm (3.8)} 
$\phi=\phi(t,x)$ satisfying for any $T>0$, 
\begin{eqnarray*}
\left\{\begin{array}{ll}
\phi \in C^0\bigl( \, [\, 0,T\, ]\, ;L^2 \bigr), \\[5pt]
\partial _x \phi 
              \in L^{\infty} \bigl( \, 0,T \, ;L^{p+1} \bigr), \\[5pt]
\partial_x \left( \, \bigl| \, 
\partial_x \bigl( \, U + U^r + \phi \, \bigr)\, 
\bigr|^{p-1} \partial_x\bigl( \, U + U^r + \phi \, \bigr) \, \right)
              \in L^{2}\bigl( \, 0,T \, ;L^{2} \bigr). 
\end{array} 
\right.\, 
\end{eqnarray*} 
Furthermore, the solution satisfies the uniform boundedness 
\begin{align*}
\begin{aligned}
\sup_{t\in [\, 0,\infty ), x\in \mathbb{R}}|\, \phi (t,x) \,| 
\le  \widetilde{C} 
     \; 
     \Bigl( \, 
     = 
     \Vert \,\phi_0 \, \Vert_{L^{\infty}} 
     + 2\, |\, u_{-} \, | + 2\, |\, u_{+} \, | \, 
\Bigr), 
\end{aligned}
\end{align*}
and also satisfies for any $T>0$, 
\begin{align*}
\begin{aligned}
&\,\Vert \, 
 \phi(t) \, \Vert_{L^2}^2 
 + \int ^{t }_{0 } \int ^{\infty}_{-\infty} 
 {\phi}^2 \, 
 \bigl( \, \partial _x U + \partial _x U^r \, \bigr) \, 
 \mathrm{d}x \mathrm{d}\tau  \\
&\qquad \qquad \: 
 + \int ^{t }_{0 }\int ^{\infty }_{-\infty } 
   \bigl| \, \partial_x u \, \bigr|^{p+1} \, 
  \mathrm{d}x \mathrm{d}\tau 
\leq \tilde{C}_{I} 
\qquad \bigl( \, t \in [\,0,T\,] \, \bigr),  
\end{aligned}
\end{align*}
where 
\begin{align*}
\begin{aligned}
&\tilde{C}_{I}
=
 C\biggl( \, \, 
 T, \: \widetilde{C} , \: 
 \| \, \phi_0 \, \| _{L^2} 
 \, \biggr), 
\end{aligned}
\end{align*}
and 
\begin{align*}
\begin{aligned}
&\,\Vert \, \partial_x 
u(t) \, \Vert_{L^{p+1}}^{p+1} 
 + \int ^{t }_{0 }\int ^{\infty }_{-\infty } 
      \bigl| \, \partial_x u \, 
      \bigr|^{2(p-1)} 
      \bigl( \, \partial_x^2 u \, 
      \bigr)^2 \, 
   \mathrm{d}x \mathrm{d}\tau \\ 
&\qquad \qquad \quad \quad  
 + \int ^{t }_{0 }\int ^{\infty }_{-\infty } 
      \bigl| \, \partial_x u \, 
      \bigr|^{p+2}  \, 
   \mathrm{d}x \mathrm{d}\tau 
\leq \tilde{C}_{I\hspace{-.1em}I} 
\qquad \bigl( \, t \in [\,0,T\,] \, \bigr), 
\end{aligned}
\end{align*}
where 
\begin{align*}
\begin{aligned}
 \tilde{C}_{I\hspace{-.1em}I}
=
 C\biggl( \, \, 
 T, \: \widetilde{C} , \: 
 \| \, \phi_0 \, \| _{L^2}, \: 
\| \, 
\partial_x u_0 
\, \| _{L^{p+1}} 
\, \biggr).  
\end{aligned}
\end{align*}
}
\noindent
Indeed, for the initial data which satisfies 
$$
\phi_0 \in L^2, 
\qquad \partial_x \phi_0 \in L^2 \cap L^{p+1}, 
$$
we can take a subsequence of $\bigl\{ u_{\epsilon } \bigr\}$ 
(write $\bigl\{ u_{\epsilon } \bigr\}$ again for simplicity) 
and a limit function $\phi$ 
(correspondingly $u := \phi + U + U^r$) 
by Lemmas 3.2, 3.3 and 3.4, 
such that 
\begin{eqnarray*}
 \left\{\begin{array}{ll}
 \displaystyle{\slim\displaylimits_{\epsilon \rightarrow 0}} \, 
 \phi_{\epsilon} = \phi 
 \in C^{0}\bigl(\, [\, 0,T\, ] \, ;L^{2} \bigr), \\[10pt]
 \displaystyle{\w*lim\displaylimits_{\epsilon \rightarrow 0}} \, 
 \partial_x \phi_{\epsilon} = \partial_x \phi
 \in L^{\infty}\bigl( \, 0,T \, ;L^{p+1} \bigr), \\[10pt]
 \displaystyle{\w*lim\displaylimits_{\epsilon \rightarrow 0}} \, 
 \partial_x u_{\epsilon} = \partial_x u 
 \in L^{\infty}\bigl( \, 0,T \, ;L^{p+1} \bigr), \\[10pt]
 \displaystyle{\wlim\displaylimits_{\epsilon \rightarrow 0}} \, 
 \partial_x \left( \, 
 \left| \, \partial_x u_{\epsilon} \, \right|^{p-1} \partial_x u_{\epsilon} \, 
 \right)
 =\partial_x \left( \, 
  \left| \, \partial_xu \, \right|^{p-1} \partial_xu \, 
  \right)
 \in L^{2}\bigl( \, 0,T \, ;L^{2} \bigr), \\[10pt]
 \displaystyle{\wlim\displaylimits_{\epsilon \rightarrow 0}} \, 
 \partial_x 
      \left( 
      \left( 
      \bigl( \, \partial_x u_{\epsilon} \, 
      \bigr)^2 + \epsilon \, 
      \right)^{\frac{p-1}{2}} \partial_x u_{\epsilon} \, 
      \right) 
 =\partial_x \left( \, 
  \left| \, \partial_xu \, \right|^{p-1} \partial_xu \, 
  \right)
 \in L^{2}\bigl( \, 0,T \, ;L^{2} \bigr). 
 \end{array}
 \right.\,
\end{eqnarray*}
We can also see that the limit function $\phi$ gives the 
unique global solution of (3.8) 
and the results in Theorem 3.3 hold. 
In particular, we note that the energy estimates in Theorem 3.3 
are independent of $\Vert \, \partial_x \phi_0 \, \Vert_{L^2}$. 
For the initial data which satisfies 
$$
\phi_0 \in L^2, 
\qquad \partial_x \phi_0 \in L^{p+1}, 
$$
we take again the approximate sequence $\bigl\{ \phi_{0}^{\delta} \bigr\}$ 
which satisfies 
$$
\phi_{0}^{\delta} \in L^2, 
\qquad \partial_x \phi_{0}^{\delta} \in L^2 \cap L^{p+1} 
$$
and 
\begin{eqnarray*}
 \left\{\begin{array}{ll}
 \displaystyle{\slim\displaylimits_{\delta \rightarrow 0}} \, 
 \phi_{0}^{\delta} = \phi_{0} 
 \in L^{2}, \\[10pt]
 \displaystyle{\slim\displaylimits_{\delta \rightarrow 0}} \, 
 \partial_x \phi_{0}^{\delta} = \partial_x \phi_{0} 
 \in L^{p+1}.
 \end{array}
 \right.\,
\end{eqnarray*}
We may take the limit $\delta \rightarrow 0$ to get Theorem 3.3. 

Since the energy estimates in Theorem 3.3 
depend on $T$, 
we can not prove the asymptotics 
$$
\Vert \, \phi (t) \, \Vert_{L^{\infty}} \rightarrow 0 
\qquad (t \rightarrow \infty). 
$$
In order to show the desired asymptotics, 
we show the following a priori estimates 
which are independent of $T$ in the next sections. 

\medskip

\noindent
{\bf Proposition 3.1} (uniform estimates 
I){\bf .}\quad {\it
For any initial data 
$\phi_0 \in L^2$ and 
$\partial _x \phi_0 \in L^{p+1}$, 
there exists a positive constant 
$$
C_{p}(\phi_0
)=
C\bigl( \, \|\phi_0\|_{L^2}
 \, \bigr)
$$ 
such that 
the unique global solution in time $\phi$ 
to the Cauchy problem {\rm (3.8)} constructed in Theorem 3.3 
satisfies 
\begin{align}
\begin{aligned}
&\| \, \phi 
(t)\, \| _{L^2}^2
+\int _0^t G(\tau) \, \mathrm{d}\tau \\
&+\int _0^t \int _{-\infty}^{\infty} 
\bigl( \partial _x \phi \bigr)^2 
\left( 
\bigl| \partial _x \phi \bigr|^{p-1} 
+ 
 \bigl| \partial _x U \bigr|^{p-1} + \bigl|  \partial _x U^r  \bigr|^{p-1}  
\right) \, \mathrm{d}x \mathrm{d}\tau 
\leq C_{p}(\phi_0)
\end{aligned}
\end{align}
for $t \geq 0$, 
where 
\begin{align*}
\begin{aligned}
G(t):=\left( \, \int_{ \tilde{U} \geq 0} 
     \phi^2 \, \partial _x \tilde{U}\, \mathrm{d}x\right) (t)   
     &+\left( \, \int_{\tilde{U} + \phi \geq 0, \tilde{U} < 0} 
     \bigl( \, \tilde{U}+\phi \, \bigr)^2 \partial _x \tilde{U}
     \, \mathrm{d}x\right) (t)\\
    &+\left( \, \int_{\tilde{U} + \phi < 0, \tilde{U} \geq 0} 
     \bigl( \, \tilde{U}+|\, \phi \, |\,  \bigr)^2 \partial _x \tilde{U}
     \, \mathrm{d}x\right) (t).
\end{aligned}
\end{align*}
}

Furthermore, we have the $L^{p+1}$-energy estimate for $\partial _x u$ as follows.  

\medskip

\noindent
{\bf Proposition 3.2} (uniform estimates 
I\hspace{-.1em}I){\bf .}\quad {\it
For any initial data 
$\phi_0 \in L^2$ and 
$\partial _x \phi_0 \in L^{p+1}$, 
there exists a positive constant 
$$
C_{p}(\phi_0, \partial _x u_0)
=
C\bigl( \, \|\phi_0\|_{L^2}, \: 
\|\partial _x u_0\|_{L^{p+1}} \, \bigr)
$$ 
such that for $t \geq 0$, 
\begin{align}
\begin{aligned}
\| \, \partial _x u 
(t)\, \| _{L^{p+1}}^{p+1} 
&+\int _0^t \int _{-\infty}^{\infty} 
 \bigl| \, \partial _x u \, \bigr|^{2(p-1)} 
 \left( \, \partial _x^2 u \, \right)^2 
 \, \mathrm{d}x \mathrm{d}\tau \\
&+\int _0^t 
 \| \, 
 \partial _x u 
 (\tau) 
 \, \| _{L^{p+2}\left( \left\{x \in \mathbb{R}\, |\, u>0 \right\} \right)}^{p+2} 
 \, \mathrm{d}\tau 
\leq C_{p}(\phi_0, \partial _x u_0). 
\end{aligned}
\end{align}
}

\medskip

{\bf Proof of Lemma 3.2.} 
We first note 
\begin{eqnarray}
\left\{\begin{array}{ll}
\displaystyle{ \| \, \phi_0^{\epsilon} \, \| _{L^q} } 
\leq \| \, \phi_0 \, \| _{L^q} 
\qquad \qquad \qquad \qquad \qquad \quad \: \: \, \bigl(1 \leq q \leq \infty \bigr), \\[5pt]
\displaystyle{ \| \, \partial_x u_0^{\epsilon} \, \| _{L^{r+1}} } 
\leq \| \, \partial_x u_0 \, \| _{L^{r+1}} 
\qquad \qquad \qquad \qquad \quad \; \: \: \, \, \, (r \geq 1), \\[5pt]
\displaystyle{ 
\sup_{- \widetilde{C}  \leq u \leq  \widetilde{C} } 
\bigl| \, D^k f^{\epsilon} (u) \, \bigr| } 
\leq \displaystyle{ 
     \sup_{0 \leq u \leq \widetilde{C} +1 } 
     \bigl| \, D^k f (u) \, \bigr| } 
\quad \quad \; \, (k=0,1,2), 
\end{array} 
\right.\,
\end{eqnarray}
where the positive constant $\widetilde{C}$ is defined in Lemma 3.1. 
By using Lemma 1.2.1, we can get 

\medskip

\noindent
{\bf Lemma 3.5.}\quad{\it

\noindent
{\rm (1)}\ $U^{r, \epsilon}(t,x)$ is 
the unique $C^\infty$-global solution in space-time 
of the Cauchy problem
$$
\left\{
\begin{array}{l} 
\partial _t U^{r, \epsilon} 
+ \partial _x \bigl( f^{\epsilon}(U^{r, \epsilon} ) \bigr) = 0 
\qquad \qquad \qquad \qquad \qquad \qquad \qquad  (t>0, x\in \mathbb{R}),\\[5pt]
U^{r, \epsilon} (0,x) = \left( \bigl( f^{\epsilon} \bigr)' \right)^{-1} 
                        \left( {\displaystyle{ 
                        \frac{\bigl( f^{\epsilon} \bigr)'(w_+)}{2} 
                        + \frac{\bigl( f^{\epsilon} \bigr)'(w_+)}{2} \tanh x 
                        }}\right) 
\; \: \: \; \; \quad \quad( x\in \mathbb{R}),\\[16pt]
\displaystyle{\lim_{x\to + \infty}} U^{r, \epsilon}(t,x) 
= \bigl( f^{\epsilon} \bigr)'(w_+)
\qquad \quad \; 
\: \: \: \quad \quad \qquad \qquad \qquad \qquad \qquad \bigl(t\ge 0 \bigr),\\[5pt]
\displaystyle{\lim_{x\to - \infty}} U^{r, \epsilon}(t,x) 
= \bigl( f^{\epsilon} \bigr)'(0)
\qquad \quad \; 
\: \; \quad \quad \quad \qquad \qquad \qquad \qquad \qquad \bigl(t\ge 0 \bigr).
\end{array}
\right.\,     
$$
{\rm (2)}\ \ 
$\bigl( f^{\epsilon} \bigr)'(0) \le U^{r, \epsilon}(t,x) \le \bigl( f^{\epsilon} \bigr)'(w_+)$ 
and\ \ $\partial_x U^{r, \epsilon}(t,x) \ge 0$  
\qquad  $(t>0, x\in \mathbb{R})$.

\smallskip

\noindent
{\rm (3)}\ For any $1\leq q \leq \infty$, there exists a positive 
constant $C_q$ such that
             \begin{eqnarray*}
                 \begin{array}{l}
                    \parallel \partial_x U^{r, \epsilon}(t)\parallel_{L^q} \leq 
                    C(q, \, \widetilde{C}) \, (1+t)^{-1+\frac{1}{q}} 
                    \; \quad \bigl(t\ge 0 \bigr),\\[5pt]
                    \parallel \partial_x^2 U^{r, \epsilon}(t) \parallel_{L^q} \leq 
                    C(q, \, \widetilde{C}) \, (1+t)^{-1}
                    \, \, \quad \quad \bigl(t\ge 0 \bigr).
                    \end{array}       
              \end{eqnarray*}
              
\smallskip

\noindent
{\rm (4)}\ It holds that 
$$
\: \displaystyle{\lim_{t\to \infty} 
\sup_{x\in \mathbb{R}} \, 
\Biggl| \,U^{r, \epsilon}(t,x) 
- \left( \bigl( f^{\epsilon} \bigr)' \right)^{-1} 
  \biggl( \, w^r \left( \, 
  \frac{x}{t} \: ;\: f^{\epsilon}(0), \, f^{\epsilon}(w_+) 
  \, \right) \biggr) 
\, \Biggr| = 0},
$$
where $w^r$ is the solution of \rm{(1.2.1)} with $w_- = 0$. 
}

\bigskip

On the other hand, we easily have 
\begin{align*}
\begin{aligned}
|\, U^{\epsilon} (t,x) \,|
&= \int_{-\infty}^{\infty} 
   \rho_{\epsilon} (y) \, 
   |\, U (t,x-y) \,|\, \mathrm{d}y \\
&\leq \int_{-\infty}^{\infty} 
      \rho_{\epsilon} (y) \, \mathrm{d}y \, 
      \| \, U 
      (t) \, \| _{L^{\infty}} 
      \leq | \, u_- \,|, 
\end{aligned}
\end{align*}
and 
\begin{align}
\begin{aligned}
\partial _x U^{\epsilon} (t,x) 
&= \int_{-\infty}^{\infty} 
   \rho_{\epsilon} (y) \, 
   \partial _x U (t,x-y) \, \mathrm{d}y \\
&\leq \int_{-\infty}^{\infty} 
      \rho_{\epsilon} (y) \, \mathrm{d}y \, 
      \| \, \partial _x U
      (t) \, \| _{L^{\infty}}, 
\end{aligned}
\end{align}
then we easily have by Lemma 2.3, 
\begin{align}
\begin{aligned}
\| \, \partial _x U^{\epsilon} 
(t) \, \| _{L^{\infty}} 
&\leq \| \, \partial _x U 
(t) \, \| _{L^{\infty}} \\
&\leq \left(2A \right)^{\frac{1}{p-1}} \, (1+t)^{- \frac{1}{p+1}} 
\leq \left(2A \right)^{\frac{1}{p-1}}. 
\end{aligned}
\end{align}
We also have for $q \geq 1$, 
\begin{align}
\begin{aligned}
\int_{-\infty}^{\infty} 
\bigl( \, 
\partial _x U^{\epsilon} \, 
\bigr)^{q} \mathrm{d}x 
&\leq \| \, \partial _x U^{\epsilon} 
      (t) \, \| _{L^{\infty}}^{q-1} 
      \int_{-\infty}^{\infty} 
      \rho_{\epsilon} (y) \, 
      \Bigl[ \, 
      U (t,x-y)\, 
      \Bigr]_{x\rightarrow -\infty }^{x\rightarrow +\infty }\, 
      \mathrm{d}y \\
&\leq \left(2A \right)^{\frac{q-1}{p-1}} |\, u_-\,| \, 
      (1+t)^{- \frac{q-1}{p+1}}, 
\end{aligned}
\end{align}
then we get 
\begin{align}
\begin{aligned}
\| \, \partial _x U^{\epsilon} 
(t)  \, \| _{L^{q}}
&\leq \left(2A \right)^{\frac{q-1}{(p-1)q}} |\, u_-\,|^{\frac{1}{q}} \, 
      (1+t)^{- \frac{q-1}{(p+1)q}} \\
&\sim \| \, \partial _x U 
      (t) \, \| _{L^{q}}. 
\end{aligned}
\end{align}
Multiplying the equation in (3.9) by 
$\phi_{\epsilon}$, 
we obtain the 
divergence form 
\begin{align}
\begin{aligned}
&\partial_t\left(\frac{1}{2} \left|\, \phi_{\epsilon} \, \right|^2 \right) 
+\displaystyle{ 
  \partial_x \biggl( \, \phi_{\epsilon} \, 
  \Bigl( \, 
  f^{\epsilon} (U^{\epsilon}+U^{r, \epsilon}+\phi_{\epsilon}) 
  - f^{\epsilon}(U^{\epsilon}+U^{r, \epsilon}) \, 
  \Bigr) \, \biggr) } \\
&+\partial _x \left( 
  - \int _{U^{\epsilon} + U^{r, \epsilon} } 
         ^{U^{\epsilon} + U^{r, \epsilon} + \phi_{\epsilon}} 
    f^{\epsilon}(s)
    \, \mathrm{d}s 
  + f^{\epsilon}\bigl( U^{\epsilon} + U^{r, \epsilon} \bigr)\, \phi_{\epsilon} 
  \, \right) \\
&+\partial _x 
  \Biggl( \, 
  -\mu \, \phi_{\epsilon} 
  \left( 
      \left( 
      \bigl( \, \partial_x u_{\epsilon} \, 
      \bigr)^2 + \epsilon \, 
      \right)^{\frac{p-1}{2}} \partial_x u_{\epsilon} 
      - \bigl( \, \partial_x U^{\epsilon} \, \bigr)^p\, 
      \right) \, 
  \Biggr) \\
&+\left( 
  f^{\epsilon}(U^{\epsilon} + U^{r, \epsilon}+\phi_{\epsilon}) 
  - f^{\epsilon}(U^{\epsilon} + U^{r, \epsilon}) 
  - 
    \bigl( f^{\epsilon} \bigr)' (U^{\epsilon} + U^{r, \epsilon}) \, \phi_{\epsilon} 
  \right) \\
&\qquad \qquad \qquad \qquad \qquad \qquad \qquad \qquad \quad \; \: \, \times 
  \bigl( \, \partial _x U^{\epsilon} + \partial _x U^{r, \epsilon} \, \bigr) \\ 
&+\mu \, \bigl( \, \partial _x \phi_{\epsilon} \, \bigr)
  \left( 
      \left( 
      \bigl( \, \partial_x u_{\epsilon} \, 
      \bigr)^2 + \epsilon \, 
      \right)^{\frac{p-1}{2}} \partial_x u_{\epsilon} \, 
      \right) \\
&= -\phi_{\epsilon} \, 
   \partial _x 
   \left( \, 
   f^{\epsilon}(U^{\epsilon} + U^{r, \epsilon}) 
   - f^{\epsilon}( U^{r, \epsilon}) \, 
   \right) 
   +\mu \, \bigl( \, \partial _x \phi_{\epsilon} \, \bigr)
    \left( \bigl( \, \partial_x U^{\epsilon} \, \bigr)^p\, \right). 
\end{aligned}
\end{align}
Integrating (3.18) with respect to $x$ and $t$, we have
\begin{align}
\begin{aligned}
&\frac{1}{2} 
\,\Vert \, \phi_{\epsilon}
(t) \, \Vert_{L^2}^2 \\
&+\int ^{t }_{0 } \int ^{\infty }_{-\infty } 
  \left( 
  f^{\epsilon}(U^{\epsilon} + U^{r, \epsilon}+\phi_{\epsilon}) 
  - f^{\epsilon}(U^{\epsilon} + U^{r, \epsilon}) 
  - 
    \bigl( f^{\epsilon} \bigr)' 
    (U^{\epsilon} + U^{r, \epsilon}) \, \phi_{\epsilon} 
  \right) \\
& \qquad \qquad \qquad \qquad \qquad \qquad \qquad \qquad \qquad \qquad \times 
  \bigl( \, \partial _x U^{\epsilon} + \partial _x U^{r, \epsilon} \, \bigr) 
  \, \mathrm{d}x \mathrm{d}\tau \\
&+\mu \int ^{t }_{0 }\int ^{\infty }_{-\infty } 
  \bigl( \, \partial _x \phi_{\epsilon} \, \bigr)
  \left( 
      \left( 
      \bigl( \, \partial_x u_{\epsilon} \, 
      \bigr)^2 + \epsilon \, 
      \right)^{\frac{p-1}{2}} \partial_x u_{\epsilon} \, 
      \right) \, 
  \mathrm{d}x \mathrm{d}\tau \\
&=\frac{1}{2} 
  \,\Vert \, \phi_0^{\epsilon} \, \Vert_{L^2}^2 
  +\int ^{t }_{0 } \int ^{\infty}_{-\infty} 
   - \phi_{\epsilon} \, 
     \left( \, 
     \bigl( f^{\epsilon} \bigr)' 
     (U^{\epsilon} + U^{r, \epsilon}) 
     -  \bigl( f^{\epsilon} \bigr)' 
        (U^{r, \epsilon}) 
     \right)
   \bigl( \, \partial_x U^{r, \epsilon} \, 
   \bigr)
   \, \mathrm{d}x \mathrm{d}\tau 
  \\
& \quad \, 
  +\int ^{t }_{0 } \int ^{\infty}_{-\infty} 
   - \phi_{\epsilon} \, 
     \left( \, 
     \bigl( f^{\epsilon} \bigr)' 
     (U^{\epsilon} + U^{r, \epsilon}) 
     \right) 
     \bigl( \, \partial_x U^{\epsilon} \, \bigr)
   \, \mathrm{d}x \mathrm{d}\tau \\ 
& \quad \, 
   + \mu \int ^{t }_{0 }\int ^{\infty }_{-\infty } 
     \bigl( \, \partial _x \phi_{\epsilon} \, \bigr)
     \left( \bigl( \, \partial_x U^{\epsilon} \, \bigr)^p\, \right) \, 
     \mathrm{d}x \mathrm{d}\tau. 
\end{aligned}
\end{align}
Since the second term on the left-hand side of (3.19) 
is equal to 
\begin{align}
\begin{aligned}
&\frac{1}{2} \, \int ^{t }_{0 } \int ^{\infty}_{-\infty} 
 \bigl( f^{\epsilon} \bigr)'' 
 \left( \, \theta \, \phi_{\epsilon} + U^{\epsilon} + U^{r, \epsilon} \, \right) 
 {\phi_{\epsilon}}^2 \, 
 \bigl( \, \partial _x U^{\epsilon} + \partial _x U^{r, \epsilon} \, \bigr) \, 
 \mathrm{d}x \mathrm{d}\tau \\
&\qquad \qquad \qquad \qquad \qquad \qquad \qquad \qquad \quad 
 (\exists 
 \theta = \theta (t,x) \in (0,1)), 
\end{aligned}
\end{align}
and since 
$\partial _x U^{\epsilon} + \partial _x U^r \geq 0$ 
by Lemma 2.2 and Lemma 2.3, 
the term is nonnegative. 
The third term on the left-hand side is also equal to 
\begin{align}
\begin{aligned}
\mu \int ^{t }_{0 }\int ^{\infty }_{-\infty } 
  \Bigl( \, \partial _x u_{\epsilon} 
            -  \bigl( \, \partial _x U^{\epsilon} + \partial _x U^{r, \epsilon} \, \bigr)
            \, 
  \Bigr)
  \left( 
      \left( 
      \bigl( \, \partial_x u_{\epsilon} \, 
      \bigr)^2 + \epsilon \, 
      \right)^{\frac{p-1}{2}} 
      \partial_x u_{\epsilon} \, 
      \right) \, 
\mathrm{d}x \mathrm{d}\tau 
\end{aligned}
\end{align}
and we can estimate it as 
\begin{align}
\begin{aligned}
&\mu \int ^{t }_{0 }\int ^{\infty }_{-\infty } 
 \left( 
 \bigl( \, \partial_x u_{\epsilon} \, 
 \bigr)^2 + \epsilon \, 
 \right)^{\frac{p-1}{2}} 
 \partial _x u_{\epsilon} \, 
 \bigl( \, \partial _x U^{\epsilon} + \partial _x U^{r, \epsilon} \, \bigr) \, 
 \mathrm{d}x \mathrm{d}\tau \\
&\leq \frac{\mu}{2} 
      \int ^{t }_{0 }\int ^{\infty }_{-\infty } 
      \left( 
      \bigl( \, \partial_x u_{\epsilon} \, 
      \bigr)^2 + \epsilon \, 
      \right)^{\frac{p-1}{2}} 
      \Bigl( 
      \bigl( \, \partial _x u_{\epsilon} \, \bigr)^2 + 
      \bigl( \, \partial _x U^{\epsilon} + \partial _x U^{r, \epsilon} \, \bigr)^2 
      \Bigr) \, 
      \mathrm{d}x \mathrm{d}\tau \\
&\leq \frac{\mu}{2} 
      \int ^{t }_{0 }\int ^{\infty }_{-\infty } 
      \left( 
      \bigl( \, \partial_x u_{\epsilon} \, 
      \bigr)^2 + \epsilon \, 
      \right)^{\frac{p-1}{2}} 
      \bigl( \, \partial _x u_{\epsilon} \, \bigr)^2  \, 
      \mathrm{d}x \mathrm{d}\tau \\
& \quad \, 
  + \int ^{t }_{0 }\int ^{\infty }_{-\infty } 
    \frac{\mu}{2} \, 2^{\frac{p-1}{2}}
    \Bigl( \, 
    \bigl| \, \partial _x u_{\epsilon} \, \bigr|^{p-1} 
    + {\epsilon}^{\frac{p-1}{2}} 
    \Bigr) 
    \bigl( \, \partial _x U^{\epsilon} + \partial _x U^{r, \epsilon} \, \bigr)^2 \, 
    \mathrm{d}x \mathrm{d}\tau \\
&\leq \frac{3}{4} \, \mu
      \int ^{t }_{0 }\int ^{\infty }_{-\infty } 
      \left( 
      \bigl( \, \partial_x u_{\epsilon} \, 
      \bigr)^2 + \epsilon \, 
      \right)^{\frac{p-1}{2}} 
      \bigl( \, \partial _x u_{\epsilon} \, \bigr)^2  \, 
      \mathrm{d}x \mathrm{d}\tau \\
& \quad \, 
  + \frac{\mu}{2} \, 2^{\frac{p-1}{2}} {\epsilon}^{\frac{p-1}{2}} 
    \int ^{t }_{0 }\int ^{\infty }_{-\infty } 
    \bigl( \, \partial _x U^{\epsilon} + \partial _x U^{r, \epsilon} \, \bigr)^2 \, 
    \mathrm{d}x \mathrm{d}\tau \\
& \quad \, 
  + 
    C_{\mu} 
    \int ^{t }_{0 }\int ^{\infty }_{-\infty } 
    \bigl( \, \partial _x U^{\epsilon} + \partial _x U^{r, \epsilon} \, \bigr)^{p+1} \, 
    \mathrm{d}x \mathrm{d}\tau \\
&\leq \frac{3}{4} \, \mu
      \int ^{t }_{0 }\int ^{\infty }_{-\infty } 
      \left( 
      \bigl( \, \partial_x u_{\epsilon} \, 
      \bigr)^2 + \epsilon \, 
      \right)^{\frac{p-1}{2}} 
      \bigl( \, \partial _x u_{\epsilon} \, \bigr)^2  \, 
      \mathrm{d}x \mathrm{d}\tau \\
& \qquad \qquad \qquad \qquad \qquad \quad 
      + C_{p,\mu,u_{\pm}}\, {\epsilon}^{\frac{p-1}{2}}T + C_{p,\mu,u_{\pm}}\, T. 
\end{aligned}
\end{align}
Substituting (3.20) and (3.21) with (3.22) into (3.19), we have 
\begin{align}
\begin{aligned}
&\frac{1}{2} 
\,\Vert \, \phi_{\epsilon}
(t) \, \Vert_{L^2}^2 
+\frac{1}{2} \, C_{p,{\epsilon}}^{-1}\, 
 \int ^{t }_{0 } \int ^{\infty}_{-\infty} 
 {\phi_{\epsilon}}^2 \, 
 \bigl( \, \partial _x U^{\epsilon} + \partial _x U^{r, \epsilon} \, \bigr) \, 
 \mathrm{d}x \mathrm{d}\tau  \\
&\qquad \qquad \quad \, \, 
 +\frac{\mu}{4} \int ^{t }_{0 }\int ^{\infty }_{-\infty } 
\bigl( \, \partial_x u_{\epsilon} \, \bigr)^2
      \left( 
      \bigl( \, \partial_x u_{\epsilon} \, 
      \bigr)^2 + \epsilon \, 
      \right)^{\frac{p-1}{2}} \, 
  \mathrm{d}x \mathrm{d}\tau \\
&\leq \frac{1}{2} 
      \,\Vert \, \phi_0^{\epsilon} \, \Vert_{L^2}^2 
      + C_{p,\mu,u_{\pm}}\, {\epsilon}^{\frac{p-1}{2}}T + C_{p,\mu,u_{\pm}}\, T \\
& \quad \, 
 + \left| \, 
   \int ^{t }_{0 } \int ^{\infty}_{-\infty} 
      - \phi_{\epsilon} \, 
        \left( \, 
        \bigl( f^{\epsilon} \bigr)'
        (U^{\epsilon} + U^{r, \epsilon}) 
        -  
           \bigl( f^{\epsilon} \bigr)'
         (U^{r, \epsilon}) 
     \right)
    \bigl( \, \partial_x U^{r, \epsilon} \, 
    \bigr)
    \, \mathrm{d}x \mathrm{d}\tau \, \right|
  \\
& \quad \, 
  +\left| \, 
  \int ^{t }_{0 } \int^{\infty}_{-\infty} 
   - \phi_{\epsilon} \, 
     \left( \, 
     \bigl( f^{\epsilon} \bigr)'
     (U^{\epsilon} + U^{r, \epsilon}) 
     \right) 
     \bigl( \, \partial_x U^{\epsilon} \, \bigr)
   \, \mathrm{d}x \mathrm{d}\tau \, \right| \\ 
& \quad \, 
   + \mu \, \left| \, \int ^{t }_{0 }\int ^{\infty }_{-\infty } 
     \bigl( \, \partial _x \phi_{\epsilon} \, \bigr)
     \left( \bigl( \, \partial_x U^{\epsilon} \, \bigr)^p\, \right) \, 
     \mathrm{d}x \mathrm{d}\tau \, \right|. 
\end{aligned}
\end{align}
Next, by using Lemma 2.2, Lemma 3.1 and (3.15), 
we estimate the second term on the right-hand side of (3.23) as 
\begin{align}
\begin{aligned}
&\left| \, 
 \int ^{t }_{0 } \int ^{\infty}_{-\infty} 
   - \phi_{\epsilon} \, 
     \left( \, 
     \bigl( f^{\epsilon} \bigr)'
     (U^{\epsilon} + U^{r, \epsilon}) 
     -  
        \bigl( f^{\epsilon} \bigr)'
        (U^{r, \epsilon}) 
     \right)
   \bigl( \, \partial_x U^{r, \epsilon} \, 
   \bigr)
   \, \mathrm{d}x \mathrm{d}\tau \, 
 \right| \\ 
&\leq \sup_{t\in [\, 0,\infty ), x\in \mathbb{R}}|\, 
      \phi_{\epsilon}(t,x)\,| 
      \displaystyle{ 
      \sup_{-\widetilde{C} \leq u \leq \widetilde{C} } 
      \left|\, 
      \bigl( f^{\epsilon} \bigr)'' (u) \, \right| } \, \\
& \qquad \qquad \qquad \qquad \qquad \times 
      \int ^{t }_{0 } \int ^{\infty}_{-\infty} 
      | \, U^{\epsilon} \, | \, \partial_x U^r \, 
      \mathrm{d}x \mathrm{d}\tau \\
&\leq C^{\dagger}T, 
\end{aligned}
\end{align}
where 
$$
C^{\dagger}
=
 C\left(\, 
 \widetilde{C}, \:  
 \displaystyle{ 
      \sup_{-\widetilde{C} \leq u \leq \widetilde{C} } 
      \left|\, 
      \bigl( f^{\epsilon} \bigr)'' (u) \, \right| } 
 \, \right).
$$ 
Similarly, we also estimate the third term on the right-hand side as 
\begin{align}
\begin{aligned}
&\left| \, 
 \int ^{t }_{0 } \int ^{\infty}_{-\infty} 
   - \phi_{\epsilon} \, 
     \left( \, 
     \bigl( f^{\epsilon} \bigr)'
     (U^{\epsilon} + U^{r, \epsilon}) 
     \right) 
     \bigl( \, \partial_x U^{\epsilon} \, \bigr)
   \, \mathrm{d}x \mathrm{d}\tau  \, 
 \right| \\ 
&\leq \sup_{t\in [\, 0,\infty ), x\in \mathbb{R}}|\,\phi_{\epsilon}(t,x)\,| 
      \displaystyle{ 
      \sup_{-\widetilde{C} \leq u \leq \widetilde{C} } 
      \left|\, 
      \bigl( f^{\epsilon} \bigr)' (u)
      \, \right| } \, \\
& \qquad \qquad \qquad \qquad \qquad \times 
      \int ^{t }_{0 } \int ^{\infty}_{-\infty} 
      | \, U^{\epsilon} + U^{r, \epsilon} \, | \, \partial_x U^{\epsilon} \, 
      \mathrm{d}x \mathrm{d}\tau \\
&\leq 
      C^{\dagger}T. 
\end{aligned}
\end{align}
Finally, we estimate the fourth term on the right-hand side as 
\begin{align}
\begin{aligned}
&\mu \, \left| \, \int ^{t }_{0 }\int ^{\infty }_{-\infty } 
     \bigl( \, \partial _x \phi_{\epsilon} \, \bigr)
     \left( \bigl( \, \partial_x U^{\epsilon} \, \bigr)^p\, \right) \, 
     \mathrm{d}x \mathrm{d}\tau \, \right| \\
&\leq \int ^{t }_{0 }\int ^{\infty }_{-\infty } 
      \left( \,  
      \frac{\mu}{8}\, \bigl| \, \partial_x u_{\epsilon} \, \bigr|^{p+1} 
      + C_{\mu} \, \bigl| \, \partial_x U^{\epsilon} \, \bigr|^{p+1} \, 
      \right) 
      \mathrm{d}x \mathrm{d}\tau \\
&\leq \frac{\mu}{8} 
      \int ^{t }_{0 }\int ^{\infty }_{-\infty } 
      \left( 
      \bigl( \, \partial_x u_{\epsilon} \, 
      \bigr)^2 + \epsilon \, 
      \right)^{\frac{p-1}{2}} 
      \bigl( \, \partial _x u_{\epsilon} \, \bigr)^2  \, 
      \mathrm{d}x \mathrm{d}\tau 
      + C_{\mu} \, \int ^{t }_{0 } 
        \| \, \partial _x U^{\epsilon} 
        (\tau )
        \, \| _{L^{p+1}}^{p+1} 
        \mathrm{d}\tau. 
\end{aligned}
\end{align}
Substituting (3.23), (3.24) and (3.25) with (3.16) into (3.22), 
we get the desired a priori energy inequality
\begin{align}
\begin{aligned}
&\frac{1}{2} 
\,\Vert \, \phi_{\epsilon}
(t) \, \Vert_{L^2}^2 
+\frac{1}{2} \, C_{p,{\epsilon}}^{-1}\, 
 \int ^{t }_{0 } \int ^{\infty}_{-\infty} 
 {\phi_{\epsilon}}^2 \, 
 \bigl( \, \partial _x U^{\epsilon} + \partial _x U^{r, \epsilon} \, \bigr) \, 
 \mathrm{d}x \mathrm{d}\tau  \\
&\qquad \qquad \quad \, \, 
 +\frac{\mu}{8} \int ^{t }_{0 }\int ^{\infty }_{-\infty } 
  \bigl( \, \partial_x u_{\epsilon} \, \bigr)^2
      \left( 
      \bigl( \, \partial_x u_{\epsilon} \, 
      \bigr)^2 + \epsilon \, 
      \right)^{\frac{p-1}{2}} \, 
  \mathrm{d}x \mathrm{d}\tau \\
&\leq \frac{1}{2} 
      \,\Vert \, \phi_0^{\epsilon} \, \Vert_{L^2}^2 
      +C^{\dagger \dagger}
       T,  
\end{aligned}
\end{align}
where 
\begin{align*}
\begin{aligned}
&C^{\dagger \dagger}
=
 C\Biggl(\,  
 \widetilde{C}, \: 
 \displaystyle{ 
      \sup_{-\widetilde{C} \leq u \leq \widetilde{C} } 
      \left|\, 
      \bigl( f^{\epsilon} \bigr)' (u) \, \right| }, \: 
 \displaystyle{ 
      \sup_{-\widetilde{C} \leq u \leq \widetilde{C} } 
      \left|\, 
      \bigl( f^{\epsilon} \bigr)'' (u) \, \right| }
 \, \Biggr).
\end{aligned}
\end{align*}
Thus, by noting (3.13), we complete the proof of Lemma 3.2. 

\bigskip

{\bf Proof of Lemma 3.3.} 
Multiplying the equation in (3.10) by 
$$-\partial_x 
      \left( 
      \left( 
      \bigl( \, \partial_x u_{\epsilon} \, 
      \bigr)^2 + \epsilon \, 
      \right)^{\frac{p-1}{2}} \partial_x u_{\epsilon} \, 
      \right), $$ 
we have the following divergence form 
\begin{align}
\begin{aligned}
&\partial_t \left( \,
 \int_{0}^{\partial_x u_{\epsilon}} 
 \bigl( \, s^2 + \epsilon \, \bigr)^{\frac{p-1}{2}} s 
 \, \mathrm{d}s \, 
 \right) \\
&+\partial _x \Biggl( 
  - \Bigl( \, 
    \bigl( \partial_x u_{\epsilon} \bigr)^2 + \epsilon \, 
    \Bigr)^{\frac{p-1}{2}} 
  \partial_x u_{\epsilon} \, 
  \biggl(\, \partial_t u_{\epsilon} 
  + \displaystyle{ 
  \partial_x 
  \Bigl( \, 
  f^{\epsilon} (u_{\epsilon}) \, 
  \Bigr) }  \, 
  \biggr)
  \, \Biggr) \\
&+\partial _x 
  \Biggl( \, 
  \bigl( f^{\epsilon} \bigr)'(u_{\epsilon}) \, 
  \int_{0}^{\partial_x u_{\epsilon}} 
  \bigl( \, s^2 + \epsilon \, \bigr)^{\frac{p-1}{2}} s 
  \, \mathrm{d}s \, 
  \Biggr) \\
&+ \mu \, 
   \Biggl( 
   \partial_x 
      \left( 
      \left( 
      \bigl( \, \partial_x u_{\epsilon} \, 
      \bigr)^2 + \epsilon \, 
      \right)^{\frac{p-1}{2}} \partial_x u_{\epsilon} \, 
      \right)
   \Biggr)^2 \\
&=-
   \bigl( f^{\epsilon} \bigr)''(u_{\epsilon}) 
   \Bigl( \, 
    \bigl( \partial_x u_{\epsilon} \bigr)^2 + \epsilon \, 
    \Bigr)^{\frac{p-1}{2}} 
   \bigl( \, \partial_x u_{\epsilon} \, \bigr)^3\\ 
&\quad \, 
 +
  \bigl( f^{\epsilon} \bigr)''(u_{\epsilon}) 
  \, \partial_x u_{\epsilon} 
  \int_{0}^{\partial_x u_{\epsilon}} 
  \bigl( \, s^2 + \epsilon \, \bigr)^{\frac{p-1}{2}} s 
  \, \mathrm{d}s. 
\end{aligned}
\end{align}
Integrating (3.28) with respect to $t$ and $x$, we have 
\begin{align}
\begin{aligned}
&\int_{-\infty}^{\infty}
 \int_{0}^{\partial_x u_{\epsilon}
 (t)} 
 \bigl( \, s^2 + \epsilon \, \bigr)^{\frac{p-1}{2}} s 
 \, \mathrm{d}s \mathrm{d}x \\
&+ \mu \, \lint_{0}^{t} \lint_{-\infty}^{\infty}
   \Biggl( 
   \partial_x 
      \left( 
      \left( 
      \bigl( \, \partial_x u_{\epsilon} \, 
      \bigr)^2 + \epsilon \, 
      \right)^{\frac{p-1}{2}} \partial_x u_{\epsilon} \, 
      \right)
   \Biggr)^2 
   \, \mathrm{d}x \mathrm{d}\tau \\
&\leq \int_{-\infty}^{\infty}
      \int_{0}^{\partial_x u^{\epsilon}_{0} } 
      \bigl( \, s^2 + \epsilon \, \bigr)^{\frac{p-1}{2}} s 
      \, \mathrm{d}s \mathrm{d}x \\
&\quad \, 
 +\int_{0}^{t} \int_{-\infty}^{\infty} 
  \left| \, 
  \bigl( f^{\epsilon} \bigr)''(u_{\epsilon}) \, 
  \right| 
   \Bigl( \, 
    \bigl( \, \partial_x u_{\epsilon} \, \bigr)^2 + \epsilon \, 
    \Bigr)^{\frac{p-1}{2}} 
   \bigl| \, \partial_x u_{\epsilon} \, \bigr|^3
  \, \mathrm{d}x \mathrm{d}\tau \\ 
&\quad \, 
 +\int_{0}^{t} \int_{-\infty}^{\infty} 
  \left| \, 
  \bigl( f^{\epsilon} \bigr)''(u_{\epsilon}) \, 
  \right| \, 
  \bigl| \, \partial_x u_{\epsilon} \, \bigr| \, 
  \biggl| \, 
  \int_{0}^{\partial_x u_{\epsilon}} 
  \bigl( \, s^2 + \epsilon \, \bigr)^{\frac{p-1}{2}} s 
  \, \mathrm{d}s \, \biggr| 
  \, \mathrm{d}x \mathrm{d}\tau \\ 
&\leq \int_{-\infty}^{\infty}
      \int_{0}^{\partial_x u^{\epsilon}_{0} } 
      \bigl( \, s^2 + \epsilon \, \bigr)^{\frac{p-1}{2}} s 
      \, \mathrm{d}s \mathrm{d}x \\
&\quad \, 
 +2 \displaystyle{ 
      \sup_{-\widetilde{C} \leq u \leq \widetilde{C} } 
      \left|\, 
      \bigl( f^{\epsilon} \bigr)'' (u) \, \right| } \, 
    \int_{0}^{t} \int_{-\infty}^{\infty} 
    \Bigl( \, 
    \bigl( \, \partial_x u_{\epsilon} \, \bigr)^2 + \epsilon \, 
    \Bigr)^{\frac{p-1}{2}} 
   \bigl| \, \partial_x u_{\epsilon} \, \bigr|^3
  \, \mathrm{d}x \mathrm{d}\tau. 
\end{aligned}
\end{align}
In order to estimate 
the second term on the right-hand side of above energy inequality, 
we use the following lemma. 
Since the proof is elementary 
we omit the proof. 

\medskip

\noindent
{\bf Lemma 3.6.}\quad {\it
There exists a positive constant $C_{p,q}$ such that 
for $v \in L^{2(q-p+1)}$ with $\partial_x v \in L^{2}$ $(p>1, \; q>p-1)$, 
it holds 
\begin{align*}
\begin{aligned}
\displaystyle{\sup_{x \in \mathbb{R}}} 
\bigl| \, v (x) \, \bigr| 
&\leq C_{p,q} 
      \left( \, 
      \int_{-\infty}^{\infty} 
      \bigl| \, v (x) \, \bigr|^{2(q-p+1)} 
      \, \mathrm{d}x \, 
      \right)^{\frac{1}{2(q+1)}} \\
&\qquad \times 
      \left( \, 
      \int_{-\infty}^{\infty} 
      \bigl| \, v (x) \, \bigr|^{2(p-1)} 
      \bigl( \, \partial_x v (x) \, \bigr)^2 
      \, \mathrm{d}x \, 
      \right)^{\frac{1}{2(q+1)}}. 
\end{aligned}
\end{align*}
}

\medskip

\noindent
By using Lemma 3.6 with $v=\partial_x u_{\epsilon}$ and $q=\frac{3}{2} p$, 
we can estimate the second term as 
\begin{align}
\begin{aligned}
&\int_{-\infty}^{\infty} 
    \Bigl( \, 
    \bigl( \partial_x u_{\epsilon} \bigr)^2 + \epsilon \, 
    \Bigr)^{\frac{p-1}{2}} 
   \bigl| \, \partial_x u_{\epsilon} \, \bigr|^3
  \, \mathrm{d}x \\
&\leq C_{p} 
      \left( \, 
      \int_{-\infty}^{\infty} 
      \bigl| \, \partial_x u_{\epsilon} \, \bigr|^{p+2} 
      \, \mathrm{d}x \, 
      \right)^{\frac{2}{3p+2}}
      \left( \, 
      \int_{-\infty}^{\infty} 
      \bigl| \, \partial_x u_{\epsilon} \, \bigr|^{2(p-1)} 
      \bigl( \, \partial_x^2 u_{\epsilon} \, \bigr)^2 
      \, \mathrm{d}x \, 
      \right)^{\frac{2}{3p+2}} \\
& \qquad \times 
      \int_{-\infty}^{\infty} 
      \Bigl( \, 
      \bigl( \, \partial_x u_{\epsilon} \, \bigr)^2 + \epsilon \, 
      \Bigr)^{\frac{p-1}{2}} 
      \bigl( \, \partial_x u_{\epsilon} \, \bigr)^2 
      \, \mathrm{d}x. 
\end{aligned}
\end{align}
Noting 
\begin{align}
\int_{-\infty}^{\infty} 
      \bigl| \, \partial_x u_{\epsilon} \, \bigr|^{p+2} 
      \, \mathrm{d}x
\leq C_{p} \int_{-\infty}^{\infty} 
      \Bigl( \, 
      \bigl( \, \partial_x u_{\epsilon} \, \bigr)^2 + \epsilon \, 
      \Bigr)^{\frac{p-1}{2}} 
      \bigl( \, \partial_x u_{\epsilon} \, \bigr)^2 
      \, \mathrm{d}x, 
\end{align}
we conclude 
\begin{align}
\begin{aligned}
&\int_{-\infty}^{\infty} 
    \Bigl( \, 
    \bigl( \partial_x u_{\epsilon} \bigr)^2 + \epsilon \, 
    \Bigr)^{\frac{p-1}{2}} 
   \bigl| \, \partial_x u_{\epsilon} \, \bigr|^3
  \, \mathrm{d}x \\
&\leq C_{p} 
      \left( \, 
      \int_{-\infty}^{\infty} 
      \bigl| \, \partial_x u_{\epsilon} \, \bigr|^{2(p-1)} 
      \bigl( \, \partial_x^2 u_{\epsilon} \, \bigr)^2 
      \, \mathrm{d}x \, 
      \right)^{\frac{1}{3p+1}} \\
& \qquad \qquad \times 
      \left( \, 
      \int_{-\infty}^{\infty} 
      \Bigl( \, 
      \bigl( \, \partial_x u_{\epsilon} \, \bigr)^2 + \epsilon \, 
      \Bigr)^{\frac{p-1}{2}} 
      \bigl( \, \partial_x u_{\epsilon} \, \bigr)^2 
      \, \mathrm{d}x 
      \right)^{\frac{3p+2}{3p+1}} \\ 
&\leq \frac{\mu}{2} \, 
      \int_{-\infty}^{\infty} 
      \bigl| \, \partial_x u_{\epsilon} \, \bigr|^{2(p-1)} 
      \bigl( \, \partial_x^2 u_{\epsilon} \, \bigr)^2 
      \, \mathrm{d}x \\
& \qquad \qquad + 
      C_{p, \mu} \, 
      \left( \, 
      \int_{-\infty}^{\infty} 
      \Bigl( \, 
      \bigl( \, \partial_x u_{\epsilon} \, \bigr)^2 + \epsilon \, 
      \Bigr)^{\frac{p-1}{2}} 
      \bigl( \, \partial_x u_{\epsilon} \, \bigr)^2 
      \, \mathrm{d}x \, 
      \right)^{\frac{3p+2}{3p}}. 
\end{aligned}
\end{align}
Substituting (3.32) and Lemma 3.2 into (3.28), we get 
\begin{align}
\begin{aligned}
&\int_{-\infty}^{\infty}
 \int_{0}^{\partial_x u_{\epsilon}
 (t)} 
 \bigl( \, s^2 + \epsilon \, \bigr)^{\frac{p-1}{2}} s 
 \, \mathrm{d}s \mathrm{d}x \\
&+ \frac{\mu}{2} \, \lint_{0}^{t} \lint_{-\infty}^{\infty}
   \Biggl( 
   \frac{p \, \bigl( \, \partial_x u_{\epsilon} \, \bigr)^2 + \epsilon }
   {\bigl( \, \partial_x u_{\epsilon} \, \bigr)^2 + \epsilon }
   \Biggr)^2 
      \left( 
      \bigl( \, \partial_x u_{\epsilon} \, \bigr)^2 
      + \epsilon \, 
      \right)^{p-1} 
      \bigl( \, \partial_x^2 u_{\epsilon} \, \bigr)^2
   \, \mathrm{d}x \mathrm{d}\tau \\
&\leq \int_{-\infty}^{\infty}
      \int_{0}^{\partial_x u^{\epsilon}_{0} } 
      \bigl( \, s^2 + \epsilon \, \bigr)^{\frac{p-1}{2}} s 
      \, \mathrm{d}s \mathrm{d}x \\ 
&\quad \, 
 +C_{p, \mu} \, C_{I} 
    \Biggl( \, 
    \displaystyle{ 
      \sup_{-\widetilde{C} - 1 \leq u \leq \widetilde{C} + 1 } 
      \left|\, 
      \bigl( f^{\epsilon} \bigr)'' (u) \, \right| } \, 
    \Biggr)^{\frac{3p+2}{3p}} \\
&\quad \quad \, \, \times 
    \displaystyle{\sup_{0 \leq t \leq T}} 
    \Biggl( \, 
    \int_{-\infty}^{\infty} 
    \Bigl( \, 
    \bigl( \, \partial_x u_{\epsilon} \, \bigr)^2 + \epsilon \, 
    \Bigr)^{\frac{p-1}{2}} 
    \bigl( \, \partial_x u_{\epsilon} \, \bigr)^2 
    \, \mathrm{d}x \, 
    \Biggr)^{\frac{2}{3p}}. 
\end{aligned}
\end{align}
Noting $\frac{2}{3p}<1$ and 
\begin{align}
\begin{aligned}
&\int_{-\infty}^{\infty} 
      \int_{0}^{\partial_x u_{\epsilon}
      (t)} 
      \bigl( \, s^2 + \epsilon \, \bigr)^{\frac{p-1}{2}} s 
      \, \mathrm{d}s \mathrm{d}x \\[5pt]
&\sim \left|\left| \, 
      \partial_x u_{\epsilon}(t) \, 
      \right| \right|_{L^{p+1}}^{p+1} 
      + {\epsilon}^{\frac{p-1}{2}}
        \left|\left| \, 
        \partial_x u_{\epsilon}(t) \, 
        \right| \right|_{L^{2}}^{2} \\[5pt]
&\sim \int_{-\infty}^{\infty} 
      \Bigl( \, 
      \bigl( \, \partial_x u_{\epsilon} \, \bigr)^2 + \epsilon \, 
      \Bigr)^{\frac{p-1}{2}} 
      \bigl| \, \partial_x u_{\epsilon} \, \bigr|^3
      \, \mathrm{d}x, 
\end{aligned}
\end{align}
we can complete the proof of Lemma 3.2. 

\bigskip 

\noindent
\section{Uniform estimates 
I}
In this section, we show the basic uniform energy estimates with $ p > 1 $ 
which is not depending on $T$, 
that is, Proposition 3.1. 
In what follows, 
we show Proposition 3.1 (also Proposition 3.2 in the next section) provided 
the solution is sufficiently smooth for simplicity 
so that we can clearly present the essential process 
to get the uniform estimates. 
In order to justify the estimates for the solution obtained in Thoerem 3.2, 
we may take $\epsilon$-regularization again as in Section 3, 
and take the limit $\epsilon \searrow 0$. 
Since the process is standard, we omit the details here. 
We first note the uniform boundedness of $\phi$ 
which is proved in Theorem 3.3, that is, 
\begin{align}
\begin{aligned}
\sup_{t\in [\, 0,\infty ), x\in \mathbb{R}}|\, \phi (t,x) \,| 
\le  \widetilde{C}. 
\end{aligned}
\end{align}

Now let us rewrite the basic $L^2$-energy inequality, 
that is Proposition 3.1 (uniform estimates 
I): 
\begin{align}
\begin{aligned}
&\| \, \phi 
(t)\, \| _{L^2}^2
+\int _0^t G(\tau) \, \mathrm{d}\tau \\
&+\int _0^t \int _{-\infty}^{\infty} 
\bigl( \, \partial _x \phi \, \bigr)^2 
\left( 
\bigl| \partial _x \phi \bigr|^{p-1} 
+ \bigl| \partial _x U \bigr|^{p-1} + \bigl| \partial _x U^r \bigr|^{p-1}  
\right) \, \mathrm{d}x \mathrm{d}\tau 
\leq 
     C_{p}(\phi_0
)
\end{aligned}
\end{align}
for $t \geq 0$, where $G(t)$ is defined as in Proposition 3.1. 
The proof of 
(4.2) is given by the following two lemmas. 

\medskip

\noindent
{\bf Lemma 4.1.}\quad {\it
It holds that for $t \geq 0$, 
\begin{align*}
\begin{aligned}
&\| \, \phi 
(t)\, \| _{L^2}^2
+\int _0^t G(\tau) \, \mathrm{d}\tau \\
& \: \; \; \; +\int _0^t \int _{-\infty}^{\infty} 
               \bigl( \, \partial _x \phi \, \bigr)^2 
               \left( 
               \bigl| \partial _x \phi \bigr|^{p-1} 
               + \bigl| \partial _x U \bigr|^{p-1} 
               + \bigl| \partial _x U^r \bigr|^{p-1}  
               \right) 
               \, \mathrm{d}x \mathrm{d}\tau \\
& \leq C_p \| \, \phi _0 \, \| _{L^2}^2 
       + C_p \int _0^t 
         \left(\, \| \, \phi 
         (\tau )\, \| _{L^2}^2 +1 \,  \right) 
         \left|\, \int _{-\infty}^{\infty}
         \left|\, \widetilde{F_{p} }(U,U^r) \, \right| \, \mathrm{d}x \, 
         \right|^{\frac{3p+1}{3p}}(\tau) 
         \, \mathrm{d}\tau \\ 
& \quad + C_p \int _0^t 
        \left|\, \int _{-\infty}^{\infty} 
        \bigl( \, \partial_x U + \partial_x U^{r} \, \bigr)^{p-1} 
        \left( \, \partial_x U^r \, \right)^{2}
        \, \mathrm{d}x \, 
        \right| (\tau) \, \mathrm{d}\tau . 
\end{aligned}
\end{align*}
} 

\medskip

\noindent
{\bf Lemma 4.2.}\quad {\it
It holds that 
\begin{align*}
\begin{aligned}
\int _0^{\infty} 
\left|\, \int _{-\infty}^{\infty}
\left|\, \widetilde{F_{p} }(U,U^r) \, \right| \, \mathrm{d}x \, 
\right|^{\frac{3p+1}{3p}}(t) 
\, \mathrm{d}t 
< \infty, 
\end{aligned}
\end{align*}
\begin{align*}
\begin{aligned}
\int _0^{\infty} 
\left|\, \int _{-\infty}^{\infty} 
\bigl( \, \partial_x U + \partial_x U^{r} \, \bigr)^{p-1} 
\left( \, \partial_x U^r \, \right)^{2}
\, \mathrm{d}x \, 
\right| (t) \, \mathrm{d}t 
< \infty. 
\end{aligned}
\end{align*}
}

Once Lemma 4.1 and Lemma 4.2 are proved, by 
Gronwall's inequality, 
we have the uniform boundedness
\begin{align*}
\begin{aligned}
\| \, \phi 
(t)\, \| _{L^2}^2 
\leq 
&C_p \bigl( \, \| \, \phi _0 \, \| _{L^2}^2 + 1 \, \bigr) \\
&\times 
     \exp 
     \left\{ \, 
     \int _0^{\infty} 
     \left|\, \int _{-\infty}^{\infty}
     \left|\, \widetilde{F_{p} }(U,U^r) \, \right| \, \mathrm{d}x \, 
     \right|^{\frac{3p+1}{3p}}
     \, \mathrm{d}t \, 
     \right\}
     < \infty 
\end{aligned}
\end{align*}
which easily implies (4.2), that is, Proposition 3.1. 

\medskip

{\bf Proof of Lemma 4.1.}
For $ p > 1$, multiplying the equation in (3.7) by 
$\phi$, 
we obtain the 
divergence form 
\begin{align}
\begin{aligned}
&\partial_t\left(\frac{1}{2} \left|\, \phi \, \right|^2 \right) \\
&+\partial _x \left( 
  \phi \, \bigl( \, f(\tilde{U}+\phi) - f(\tilde{U}) \, \bigr) 
  \right) \\
&+\partial _x \left( 
  - \int _{\tilde{U} }^{\tilde{U} + \phi} f(s)
    \, \mathrm{d}s 
  + f(\tilde{U} )\, \phi 
  \, \right) \\
&+\partial _x \biggl( \, 
  -\mu \, \phi \left( \, 
    \bigl| \, \partial_x \tilde{U} + \partial_x \phi \, \bigr|^{p-1} 
    \bigl( \, \partial_x \tilde{U} + \partial_x \phi \, \bigr) - 
    \bigl| \, \partial_x \tilde{U} \, \bigr|^{p-1} \partial_x \tilde{U}  \, 
    \right)
  \, \biggr) \\
&+\left( f(\tilde{U}+\phi) - f(\tilde{U}) - f'(\tilde{U}) \, \phi \right) 
  \partial _x \tilde{U}  \\ 
&+\mu \, \bigl( \, \partial _x \phi \, \bigr)
  \left( \, 
    \bigl| \, \partial_x \tilde{U} + \partial_x \phi \, \bigr|^{p-1} 
    \bigl( \, \partial_x \tilde{U} + \partial_x \phi \, \bigr) - 
    \bigl| \, \partial_x \tilde{U} \, \bigr|^{p-1} \partial_x \tilde{U}  \, 
    \right) 
= \phi \, F_{p}(U,U^r). 
\end{aligned}
\end{align}
Integrating (4.3) with respect to $x$ and $t$, we have
\begin{align}
\begin{aligned}
&\frac{1}{2} 
\,\Vert \, \phi
(t) \, \Vert_{L^2}^2 
+\int ^{t }_{0 } \int ^{\infty }_{-\infty } 
\left( f(\tilde{U}+\phi) - f(\tilde{U}) - f'(\tilde{U}) \, \phi \right) 
\partial _x \tilde{U}
\, \mathrm{d}x \mathrm{d}\tau \\
&+\mu \int ^{t }_{0 }\int ^{\infty }_{-\infty } 
\bigl( \, \partial _x \phi \, \bigr)
  \left( \, 
    \bigl| \, \partial_x \tilde{U} + \partial_x \phi \, \bigr|^{p-1} 
    \bigl( \, \partial_x \tilde{U} + \partial_x \phi \, \bigr) - 
    \bigl| \, \partial_x \tilde{U} \, \bigr|^{p-1} \partial_x \tilde{U}  \, 
    \right) 
\, \mathrm{d}x \mathrm{d}\tau \\
&=\frac{1}{2} 
\,\Vert \, \phi_0 \, \Vert_{L^2}^2 
+\int ^{t }_{0 } \int ^{\infty}_{-\infty} 
 \phi \,  \, F_{p}(U,U^r) 
 \, \mathrm{d}x \mathrm{d}\tau. 
\end{aligned}
\end{align}
To estimate the second term on the left-hand side of 
(4.4), 
noting the shape of the flux function $f$, 
we divide the integral region of $x$ depending on the signs of
$\tilde U +\phi$, $\tilde U$ and $\phi$ as 
\begin{align}
\begin{aligned}
&\int ^{\infty }_{-\infty } 
\left( f(\tilde{U}+\phi) - f(\tilde{U}) - f'(\tilde{U}) \, \phi \right) 
\partial _x \tilde{U}
\, \mathrm{d}x 
\\
\qquad &= \int ^{\infty }_{-\infty } \biggl(\,  \int _0^{\phi} 
\left( \lambda \bigl(\tilde{U}+\eta \bigr)-\lambda \bigl(\tilde{U} \bigr) \right) 
\, \mathrm{d}\eta \biggr)
\, \bigl( \, \partial _x \tilde{U} \, \bigr) \, \mathrm{d}x 
\\
\qquad &= \int_{\tilde{U}+\phi \geq 0, \tilde{U} \geq 0, \phi \geq 0}
          +\int_{\tilde{U}+\phi \geq 0, \tilde{U} \geq 0, \phi \le 0}
          +\int_{\tilde{U}+\phi \geq 0, \tilde{U} < 0}
          + \int_{\tilde{U}+\phi < 0, \tilde{U} \geq 0}
\end{aligned}
\end{align}
where we used the fact that the integral is clearly zero on the 
region $\tilde{U}+\phi \le 0$ and $\tilde{U}\le 0$.
By 
Lagrange's mean-value theorem, we easily get as 
\begin{align}
\begin{aligned}
&\left( \int ^{\infty }_{-\infty } \biggl(\,  \int _0^{\phi} 
\left( \lambda \bigl(\tilde{U}+\eta \bigr)-\lambda \bigl(\tilde{U} \bigr) \right) 
\, \mathrm{d}\eta \biggr)
\, \bigl( \, \partial _x \tilde{U} \, \bigr) \, \mathrm{d}x \right) (t)
\sim  G(t) 
\end{aligned}
\end{align}
where $G=G(t)$ is defined in Proposition 3.1 
(cf. \cite{matsumura-yoshida}, \cite{yoshida}). 
Next, we also estimate the third term on the left-hand side of (4.4) as 
\begin{align}
\begin{aligned}
&\int ^{\infty }_{-\infty } 
 \bigl( \, \partial _x \phi \, \bigr)
  \left( \, 
    \bigl| \, \partial_x \tilde{U} + \partial_x \phi \, \bigr|^{p-1} 
    \bigl( \, \partial_x \tilde{U} + \partial_x \phi \, \bigr) - 
    \bigl| \, \partial_x \tilde{U} \, \bigr|^{p-1} \partial_x \tilde{U}  \, 
    \right) 
\, \mathrm{d}x \\
&\geq \nu_{p}^{-1} \int _{-\infty}^{\infty} 
      \bigl( \, \partial _x \phi \, \bigr)^2 
      \left( 
      \bigl| \partial _x \phi \bigr|^{p-1} 
      + \bigl| \partial _x U \bigr|^{p-1} + \bigl| \partial _x U^r \bigr|^{p-1}  
      \right) \, \mathrm{d}x 
\end{aligned}
\end{align}
for some constant $\nu_{p}>0$ which is depend only on $p$.
Here, we used the following absolute inequality with $p>1$, 
for any $a,b \in \mathbb{R}$, 
\begin{align}
\begin{aligned}
&\left(\, |\, a\, |^{p-1}a - |\, b\, |^{p-1}b \, \right) 
 \left(\, a-b \, \right) \\
&\geq C_{p}^{-1} 
      \left(\, |\, a\, |^{p-1} + |\, b\, |^{p-1} \, \right) 
      \left(\, a-b \, \right)^2 \\
&\geq \widetilde{C_{p}}^{-1} 
      \left(\, |\, a\, |^{p-1} + |\, b\, |^{p-1} + |\, a-b\, |^{p-1} \, \right) 
      \left(\, a-b \, \right)^2 
\end{aligned}
\end{align}
for some $C_{p},\widetilde{C_{p}}>0$ depending only on $p$. 
Furthermore, we should note 
\begin{align}
\begin{aligned}
\left|\, 
\int _{-\infty}^{\infty} 
\phi \, F_{p}(U,U^r) \, \mathrm{d}x \right|
&\leq 
\int _{-\infty}^{\infty} 
|\, \phi \, | \left|\, \widetilde{F_{p} }(U,U^r) \, \right|
\, \mathrm{d}x \\
&+ \mu \, \int _{-\infty}^{\infty} 
  \bigl|\, \partial_x \phi \, \bigr| 
  \left( \, 
  \bigl( \, \partial_x U + \partial_x U^{r} \, \bigr)^{p} 
  - \left( \, \partial_x U \, \right)^{p} \, 
    \right) 
    \, \mathrm{d}x. 
\end{aligned}
\end{align}
Substituting (4.6), (4.7) and (4.9) into (4.4), we get the 
energy inequality 
\begin{align}
\begin{aligned}
&\frac{1}{2} \, 
\| \, \phi 
(t)\, \| _{L^2}^2
+ C_{p}^{-1} \int _0^t G(\tau) \, \mathrm{d}\tau \\
& \: \; \; \; +\mu \, \nu_{p}^{-1} \int _0^t \int _{-\infty}^{\infty} 
               \bigl( \, \partial _x \phi \, \bigr)^2 
               \left( 
               \bigl| \partial _x \phi \bigr|^{p-1} 
               + \bigl| \partial _x U \bigr|^{p-1} 
               + \bigl| \partial _x U^r \bigr|^{p-1}  
               \right) 
               \, \mathrm{d}x \mathrm{d}\tau \\
& \leq \frac{1}{2} \, \| \, \phi _0 \, \| _{L^2}^2 
       + \int _0^t \int _{-\infty}^{\infty} 
         |\, \phi \, | \left|\, \widetilde{F_{p} }(U,U^r) \, \right|
         \, \mathrm{d}x 
         \mathrm{d}\tau \\ 
& \; \quad + \mu \, \int _0^t \int _{-\infty}^{\infty} 
          \bigl|\, \partial_x \phi \, \bigr| 
          \left( \, 
          \bigl( \, \partial_x U + \partial_x U^{r} \, \bigr)^{p} 
          - \left( \, \partial_x U \, \right)^{p} \, 
          \right) 
          \, \mathrm{d}x \mathrm{d}\tau. 
\end{aligned}
\end{align}
We estimate the second term on the right-hand side of (4.10) as follows: 
\begin{align}
\begin{aligned}
&\int _{-\infty}^{\infty} 
 |\, \phi \, | \left|\, \widetilde{F_{p} }(U,U^r) \, \right|
 \, \mathrm{d}x \\
&\leq C_{p} \| \, \phi \, \| _{L^2}^{\frac{2p}{3p+1}} 
            \| \, \partial_x \phi \, \| _{L^{p+1}}^{\frac{p+1}{3p+1}}
            \int _{-\infty}^{\infty} 
            \left|\, \widetilde{F_{p} }(U,U^r) \, \right|
            \, \mathrm{d}x \\
&\leq \frac{\mu}{4\, \nu_p} \| \, \partial_x \phi \, \| _{L^{p+1}}^{p+1}
      + C_{p} \| \, \phi \, \| _{L^2}^{\frac{2}{3}} 
              \left(\, \int _{-\infty}^{\infty} 
              \left|\, \widetilde{F_{p} }(U,U^r) \, \right|
              \, \mathrm{d}x 
              \right)^{\frac{3p+1}{3p}} \\
&\leq \frac{\mu}{4\, \nu_p} \| \, \partial_x \phi \, \| _{L^{p+1}}^{p+1}
      + C_{p} \bigl( \, \| \, \phi \, \| _{L^2}^{2} + 1 \, \bigr) 
              \left|\, \int _{-\infty}^{\infty} 
              \left|\, \widetilde{F_{p} }(U,U^r) \, \right|
              \, \mathrm{d}x 
              \, \right|^{\frac{3p+1}{3p}}, 
\end{aligned}
\end{align}
where we used Young's inequality 
and the following Sobolev type inequality (cf. \cite{yoshida}): 
\begin{align}
\begin{aligned}
\| \, \phi \, \| _{L^{\infty}} 
\leq \left( \frac{3p+1}{p+1} \right)^{\frac{p+1}{3p+1}}
     \| \, \phi \, \| _{L^2}^{\frac{2p}{3p+1}} 
     \| \, \partial_x \phi \, \| _{L^{p+1}}^{\frac{p+1}{3p+1}}.
\end{aligned}
\end{align}
By the Cauchy-Schwarz inequality 
and Young's inequality, 
we also estimate the third term on the right-hand side of (4.10) as follows: 
\begin{align}
\begin{aligned}
& \mu \,  \int _{-\infty}^{\infty} 
          \bigl|\, \partial_x \phi \, \bigr| 
          \left( \, 
          \bigl( \, \partial_x U + \partial_x U^{r} \, \bigr)^{p} 
          - \left( \, \partial_x U \, \right)^{p} \, 
          \right) 
          \, \mathrm{d}x \\
& = \mu \,  p
       \int _{-\infty}^{\infty} 
       \bigl|\, \partial_x \phi \, \bigr| 
       \left( \, 
       \bigl( \, 
       \partial_x U + \theta \partial_x U^{r} 
       \, \bigr)^{p-1} 
       \partial_x U^{r} \, 
       \right) 
       \, \mathrm{d}x \\
&\qquad \qquad \qquad \qquad \qquad \quad 
 \bigl( \, \exists 
 \theta = \theta (t,x) \in (0,1) \, \bigr)\\ 
& \leq \frac{\mu}{4\, \nu_{p}} 
       \int _{-\infty}^{\infty} 
       \bigl(\, \partial_x \phi \, \bigr)^2 
       \bigl( \, 
       \partial_x U + \partial_x U^{r} 
       \, \bigr)^{p-1} 
       \, \mathrm{d}x \\
&\qquad + C_{p} 
         \left|\,  \int _{-\infty}^{\infty} 
         \bigl(\, \partial_x U^{r} \, \bigr)^2 
         \bigl( \, 
         \partial_x U + \partial_x U^{r} 
         \, \bigr)^{p-1} 
         \, \mathrm{d}x \,\right|. 
\end{aligned}
\end{align}
Thus, substituting (4.11) and (4.13) into (4.10), 
we complete the proof of Lemma 4.1. 

\bigskip

{\bf Proof of Lemma 4.2.}
Firstly, we have 
\begin{align}
\begin{aligned}
&\left|\,  \int _{-\infty}^{\infty} 
\bigl(\, \partial_x U^{r} \, \bigr)^2 
\bigl( \, 
\partial_x U + \partial_x U^{r} 
\, \bigr)^{p-1} 
\, \mathrm{d}x \,\right|\\
& \leq C_{p} 
       \| \, \partial_x \tilde{U}
       (t) \, \| _{L^{\infty}}^{p-1} 
       \| \, \partial_x U^r (t) 
       \, \| _{L^{\infty}} \, \\
& \leq C_{p} (1+t )^{-1 - \frac{p-1}{p+1}}, 
\end{aligned}
\end{align}
that is, 
\begin{align}
\begin{aligned}
\int _{-\infty}^{\infty} 
\bigl(\, \partial_x U^{r} \, \bigr)^2 
\bigl( \, 
\partial_x U + \partial_x U^{r} 
\, \bigr)^{p-1} 
\, \mathrm{d}x 
\in L_{t}^1(0,\infty) 
\end{aligned}
\end{align}  
where we used Lemma 2.2 and Lemma 2.3. 
Then, it suffices to show, 
by the definition of the remainder term $\widetilde{F_{p} }(U,U^r)$, 
that 
\begin{equation}
\int _{-\infty}^{\infty} 
\bigl| \, f'(U+U^r)-f'(U^r) \, \bigr| \, \partial _x U^r  
\, \mathrm{d}x  
\in L_t^{\frac{3p + 1}{3p}}(0,\infty),
\end{equation}
\begin{equation}
\int _{-\infty}^{\infty} \bigl| \, f'(U+U^r) \, \bigr| \, 
\partial _x U 
\, \mathrm{d}x
\in L_t^{\frac{3p + 1}{3p}}(0,\infty).
\end{equation}
To obtain (4.16) and (4.17), it is very natural to
divide the integral region $\mathbb{R}$ depending on the 
sign of $\tilde{U} = U + U^r$.
So, for any $t\ge 0$, we introduce 
$$
X:[\, 0,\infty)\ni t \longmapsto X(t)\in \mathbb{R}
$$
such that
\begin{equation}
\tilde{U}\bigl(t,X(t)\bigr) 
= U\bigl(t,X(t)\bigr)+U^r\bigl(t,X(t)\bigr)= 0
\quad \bigl(t \ge 0\bigr),
\end{equation}
that is, 
\begin{align}
\begin{aligned}
U^r\bigl(t,X(t)\bigr)
&= -U\bigl(t,X(t)\bigr) \\ 
&= \Lint^{\infty}_{X(t)}
  \frac{1}{(1+t)^{\frac{1}{p+1}}}\, 
  \Bigggl( \left( 
  A-B \left(\frac{y}{(1+t)^{\frac{1}{p+1}}} \right)^2 
  \right)\vee 0 \Bigggr)^{\frac{1}{p-1}}
  \, \mathrm{d}y. 
\end{aligned}
\end{align}
Here we note that 
$X(t)$ uniquely exists because 
$U^r$ is strictly monotonically increasing with respect to $x$ on the whole $\mathbb{R}$ 
and $U$ is also strictly monotonically increasing on 
$
-\sqrt{\frac{A}{B}}\, (1+t)^{\frac{1}{p+1}} < x < \sqrt{\frac{A}{B}}\, (1+t)^{\frac{1}{p+1}}
$.
Furthermore, note that 
$\tilde{U}(t, -\infty)= u_-<0<u_+ =\tilde{U}(t, \infty)$.
Therefore we can divide the integral region $\mathbb{R}$ into 
$\bigl(-\infty,X(t)\bigr)$ where $\tilde{U}<0$ and 
$\bigl(X(t),\infty\bigr)$ where $\tilde{U}>0$.
As a 
basic behavior of $X(t)$, we can show by 
Lemma 2.2 and Lemma 2.3 that
there exists a positive time $T_0$ such that 
for some $\delta \in \left(\, 0, \sqrt{\frac{A}{B}} \, \right)$, 
\begin{equation}
\left(\, \sqrt{\frac{A}{B}} - \delta \, \right)\, (1+t)^{\frac{1}{p+1}} 
< X(t) < \sqrt{\frac{A}{B}}\, (1+t)^{\frac{1}{p+1}}
\qquad \bigl(t \ge T_0\bigr).
\end{equation}
Indeed, 
by an easy fact 
\begin{equation}
\sup_{x\in \mathbb{R}}\, 
\left| \,u^r \left( \frac{x}{1+t}\right) 
- u^r \left(\frac{x}{t}\right) \right| 
\le 
C(1+t)^{-1}, 
\end{equation} 
and Lemma 2.2, it follows that 
\begin{equation}
\sup_{x\in \mathbb{R}}\, \left|\,U^r(t,x)- u^r \left(\frac{x}{1+t} \right) \right| 
\le 
C_\epsilon(1+t)^{-1+\epsilon} \quad \bigl(\epsilon \in (0,1)\bigr),
\end{equation} 
which implies 
\begin{align}
\begin{aligned}
&\lim_{t\to \infty}\tilde{U} 
 \bigggl(\, 
 t,\left(\, \sqrt{\frac{A}{B}} - \delta \, \right)\, (1+t)^{\frac{1}{p+1}} 
 \, \bigggr) \\
&\quad \quad =  - \lint^{\infty}_{\sqrt{\frac{A}{B}} - \delta }
                \Bigl ( \left( \, 
                A-B \, \xi^2 \, 
                \right)\vee 0 \Bigr)^{\frac{1}{p-1}}
                \, \mathrm{d}\xi <0,
\end{aligned}
\end{align}
and
\begin{equation}
\tilde{U} 
 \left(\, 
 t,\sqrt{\frac{A}{B}}\, (1+t)^{\frac{1}{p+1}} 
 \, \right) 
 = U^r 
 \left(\, 
 t,\sqrt{\frac{A}{B}}\, (1+t)^{\frac{1}{p+1}} 
 \, \right)  \, >0 
 \quad \left( \, \forall t \geq 0 \, \right).
\end{equation}
So we have (4.20) by (4.23) and (4.24).
Then, by (4.18), (4.21) and (4.22),
we have for any $\epsilon \in (0,1)$, 
there exists a positive constant $C_{\epsilon}$ such that 
\begin{equation}
\left|\, 
(\lambda)^{-1} \left(\frac{X(t)}{1+t} \right) -
\lint^{\infty}_{\frac{\mathlarger{X(t)}}{ \mathlarger{(1+t)^{\frac{1}{p+1}} }}}
                \Bigl ( \left( \, 
                A-B \, \xi^2 \, 
                \right)\vee 0 \Bigr)^{\frac{1}{p-1}}
                \, \mathrm{d}\xi 
\, \right| \le C_{\epsilon}(1+t)^{-1 + \epsilon}, 
\end{equation}
for $t \ge T_0$. 
Using (4.25), 
we can show more
precise large time behavior of $X(t)$ 
as in the following lemma.

\bigskip

\noindent
{\bf Lemma 4.3.}\quad {\it 
It holds that for each $p>1$, 
there exists a positive constant $C_{p}$ such that 
$$
\left| \, 
\sqrt{\frac{A}{B}} - \frac{X(t)}{ (1+t)^{\frac{1}{p+1}} } 
\, \right|
\leq 
     C_{p} (1+t)^{-{\frac{p-1}{p+1}}}
\qquad \bigl(t \ge T_0\bigr).
$$
}

Let us admit Lemma 4.3 for a moment 
and complete the proof of Lemma 4.2. 
We shall give the proof of Lemma 4.3 
after the proof of Lemma 4.2.
Using Lemmas 2.2, 2.3 and 4.3, we first prove (4.16).
Dividing the integral region as we mentioned above as
\begin{equation*}
\int  ^{\infty }_{-\infty }
\bigl| \, f'(U+U^r)-f'(U^r) \, \bigr| \, \partial _x U^r  
\, \mathrm{d}x  
=\int ^{X(t)}_{-\infty }
+\int ^{\infty }_{X(t)}
=: I_{11}+ I_{12},
\end{equation*}
we estimate each integral as follows: 
\begin{align}
\begin{aligned}
I_{11}(t)  & = \int ^{X(t)}_{-\infty }
               \bigl| \, f'(U+U^r)-f'(U^r) \, \bigr| \, \partial _x U^r 
               \, \mathrm{d}x \qquad \qquad\\
           & = \int ^{X(t)}_{-\infty } 
               \partial _x \bigl(f(U^r)\bigr) 
               \, \mathrm{d}x\\
           & \leq C\, \left|\, U^r\bigl(t,X(t)\bigr) \right|^2\\
           & \leq C\, \left( \, 
                  \frac{X(t)}{1+t} + C_{\epsilon}(1+t)^{-1 + \epsilon}
                  \, \right)^2\\
           & \leq C_{p}\, (1+t)^{-1-\frac{p-1}{p+1}} 
                  + C_{\epsilon}(1+t)^{-1 - \left(1-2 \epsilon \right)} 
           \qquad \bigl(\epsilon \in (0,1),\ t \ge 0\bigr),
\end{aligned}
\end{align}
\begin{align}
\begin{aligned}
I_{12}(t) & = \int ^{\infty }_{X(t)} 
              \bigl| \, f'(U+U^r)-f'(U^r) \, \bigr| \, 
              \partial _x U^r 
              \,  \mathrm{d}x\\
            & \leq  C\int ^{\infty }_{X(t)} 
              | \, U \, |\, \partial _x U^r 
              \, \mathrm{d}x\\
            & \leq  C (1+t)^{-1} \lint ^{\infty }_{X(t)} 
                    \lint ^{\infty }_{\frac{\mathlarger{x}}
                                     {\mathlarger{(1+t)}^{\frac{1}{p+1}}}} 
                    \Bigl ( \left( \, 
                    A-B \, \xi^2 \, 
                    \right)\vee 0 \Bigr)^{\frac{1}{p-1}} 
                    \, \mathrm{d}\xi \mathrm{d}x\\
            & = C (1+t)^{-1}
                \lint ^{\infty }_{\frac{\mathlarger{X(t)}}
                                 {\mathlarger{(1+t)}^{\frac{1}{p+1}}}}
                \Bigl ( \left( \, 
                A-B \, \xi^2 \, 
                \right)\vee 0 \Bigr)^{\frac{1}{p-1}} 
                \left( \, \xi \, (1+t)^{\frac{1}{p+1}} - X(t) \right) 
                \, \mathrm{d}\xi \\
            & \leq  C_{p} (1+t)^{-1 + \frac{1}{p+1}}
                    \lint ^{\sqrt{\frac{A}{B}} }_{\frac{\mathlarger{X(t)}}
                          {\mathlarger{(1+t)}^{\frac{1}{p+1}}}}
                    \xi \left( 
                    A-B \, \xi^2 
                    \right)^{\frac{1}{p-1}}
                    \, \mathrm{d}\xi \\
            & =  C_{p} (1+t)^{-1 + \frac{1}{p+1}} \cdot 
                 \frac{p-1}{2 \, B\, p} 
                 \left( \, 
                 A-B \left(\frac{X(t)}{(1+t)^{\frac{1}{p+1}}} \right)^2 \, 
                 \right)^{\frac{p}{p-1}}\\
            & \leq  C_{p} (1+t)^{-1 + \frac{1}{p+1}} 
                    \left| \, 
                    \sqrt{\frac{A}{B}} - \frac{X(t)}{(1+t)^{\frac{1}{p+1}}} 
                    \, \right|^{\frac{p}{p-1}}\\
            & \leq  C_{p} (1+t)^
                    {-1 - \frac{p-1}{p+1} }
             \qquad \bigl( t\ge T_0\bigr),
\end{aligned}
\end{align}
where we used the facts 
$\| \, \partial_x U^r 
(t) \, \|_{L^{\infty}}
\leq C(1+t)^{-1}$ in Lemma 2.2 and Lemma 4.3.
Hence, choosing $\epsilon$ suitably small in (4.26),
we can easily conclude $I_{11}, I_{12}\in L_t^{\frac{3p+1}{3p}}(0,\infty)$,
which proves (4.16). Next, we similarly show (4.17).
In this case, noting 
\begin{equation*}
\int ^{\infty }_{-\infty }f'(U+U^r)\, \partial _x U \, \mathrm{d}x =
\int ^{\infty }_{X(t)} f'(U+U^r)\, \partial _x U \, \mathrm{d}x =: I_{21},
\end{equation*}
we estimate $I_{21}$, by integration by parts, as follows:
\begin{align*}
\begin{aligned}
 I_{21}(t)=& \int ^{\infty }_{X(t)} f'(U+U^r)\, \partial _x U \, \mathrm{d}x\\
          =& - \int ^{\infty }_{X(t)} U\,  f''(U+U^r) 
             \bigl( \, \partial _x U + \partial _x U^r  \, \bigr) \, \mathrm{d}x \\
          \leq& \ C\int ^{\infty }_{X(t)} 
                  -\frac{1}{2}\, \partial _x \bigl(U^2 \bigr) 
                  \, \mathrm{d}x 
                  + C\int ^{\infty }_{X(t)} 
                    | \, U \, |\,  \partial _x U^r \, \mathrm{d}x \\
          \leq& \ C\,\left|\,U^r\bigl(t,X(t)\bigr) \right|^2 + 
                  C_p (1+t)^
                  {-1 - \frac{p-1}{p+1} } \\
          \leq& \ C_{\epsilon}
                  (1+t)^{-1 - \left(1-2 \epsilon \right)} 
                  + C_p(1+t)^{-1 - \frac{p-1}{p+1} }
\qquad \bigl(\epsilon \in (0,1),\ t\ge T_0\bigr).
\end{aligned}
\end{align*}
Hence, choosing $\epsilon $ suitably small again, we easily have
$I_{21}\in L^{\frac{3p+1}{3p}}(0,\infty)$. Thus, 
the proof of Lemma 4.2 is complete.
 
\medskip

{\bf Proof of Lemma 4.3.}\quad 
First, we note from (4.20) that 
\begin{align}
\begin{aligned}
\left(\, \sqrt{\frac{A}{B}} - \delta \, \right)\, (1+t)^{- \frac{p}{p+1}} 
< \frac{X(t)}{1+t} < \sqrt{\frac{A}{B}}\, (1+t)^{- \frac{p}{p+1}}
\qquad \bigl(t \ge T_0\bigr). 
\end{aligned}
\end{align}
Then, by using (4.25), we have for any $\epsilon \in (0,1)$, 
there exists a positive constant $C_{\epsilon}$ such that 
\begin{align}
\begin{aligned}
& (\lambda)^{-1} \left( 
  \sqrt{\frac{A}{B}}\, (1+t)^{- \frac{p}{p+1}} 
  \right) 
  +  C_{\epsilon}(1+t)^{-1 + \epsilon} \\
& \geq (\lambda)^{-1} \left(\frac{X(t)}{1+t} \right) 
       +  C_{\epsilon}(1+t)^{-1 + \epsilon}\\
& \geq B^{\frac{1}{p-1}}
       \lint^{\sqrt{\frac{A}{B}}}_{\frac{\mathlarger{X(t)}}
                                  {\mathlarger{(1+t)}^{\frac{1}{p+1}}}}
       \left( 
       \frac{A}{B} - \, \xi^2 
       \right)^{\frac{1}{p-1}}
       \, \mathrm{d}\xi \\
& \geq \left( \, A B \, \right)^{\frac{1}{2(p-1)}}
       \lint^{\sqrt{\frac{A}{B}}}_{\frac{\mathlarger{X(t)}}
                                  {\mathlarger{(1+t)}^{\frac{1}{p+1}}}}
       \left( 
       \sqrt{\frac{A}{B}} - \, \xi 
       \right)^{\frac{1}{p-1}}
       \, \mathrm{d}\xi 
\qquad \bigl(t \ge T_0\bigr). 
\end{aligned}
\end{align}
Hence, 
by taking $\epsilon = \frac{1}{p+1}$, 
it implies that for $t \ge T_0$, 
\begin{align}
\begin{aligned}
C_{p} (1+t)^{- \frac{p}{p+1}} 
\geq \left( \, A B \, \right)^{\frac{1}{2(p-1)}}
     \left( \frac{p-1}{p} \right) 
     \left| \, 
     \sqrt{\frac{A}{B}} - \frac{X(t)}{ (1+t)^{\frac{1}{p+1}} } 
     \, \right|^{\frac{p}{p-1}} 
\end{aligned}
\end{align}
which completes the proof of Lemma 4.4.

\medskip
Thus, we do complete the proof of Proposition 3.1.

\medskip

\noindent
{\bf Remark 4.1.}\quad
The unique global solution in time $u$ also satisfies 
the following regularity 
\begin{eqnarray*}
\left\{\begin{array}{ll}
\partial _t u 
\in L^{\infty} \bigl( \, 0,T \, ;L^{p+1} \bigr)
    \cap \Bigl( \, L^{p+1}\bigl(\, 0,T \, ;L^{p+1}  \bigr)
    \oplus  L^{2}\bigl( \, 0,T \, ;L^{2} \bigr) \, \Bigr),\\[5pt] 
\lambda (u) \, \partial _x u 
\in L^{\infty} \bigl( \, 0,T \, ;L^{p+1} \bigr).
\end{array} 
\right.\,
\end{eqnarray*}
. 
\bigskip 

\noindent
\section{Uniform estimates 
I\hspace{-.1em}I}
In this section, 
In order to complete the uniform estimates for the asymptotics not depending on $T$, 
we show Proposition 3.2. 
To do that, we assume that the solution to our Cauchy problem (3.8) 
satisfies the same regularity as in Section 4. 
What we should prove is 
the following energy inequality: 

\begin{align}
\begin{aligned}
\| \, \partial _x u 
(t)\, \| _{L^{p+1}}^{p+1} 
&+\int _0^t \int _{-\infty}^{\infty} 
 \bigl| \, \partial _x u \, \bigr|^{2(p-1)} 
 \left( \, \partial _x^2 u \, \right)^2 
 \, \mathrm{d}x \mathrm{d}\tau \\
&+\int _0^t 
 \| \, 
 \partial _x u 
 (\tau) 
 \, \| _{L^{p+2}\left( \left\{x \in \mathbb{R}\, |\, u>0 \right\} \right)}^{p+2} 
 \, \mathrm{d}\tau 
\leq C_{p}(\phi_0, \partial _x u_0)\quad \bigl( t \geq 0 \bigr).
\end{aligned}
\end{align}
In order to obtain (5.1), 
we multiple the equation in the original Cauchy problem (1.1) 
(not the reformulated problem, that is (3.7) or (3.8)) by 
$$
- \partial_x \left( \, 
\left| \, \partial_xu \, \right|^{q-1} \partial_xu \, 
\right)
$$ 
with $q>1$ 
and obtain the divergence form 
\begin{align}
\begin{aligned}
&\partial_t\left(\frac{1}{q+1} \left|\, \partial _x u \, \right|^{q+1} \right) 
+\partial _x \left( 
 - \, \left|\, \partial _x u \, \right|^{q-1} 
   \partial _x u \cdot \partial _t u \, \right) \\
&\qquad \; \, + \partial _x \left( - \, \frac{q}{q+1} \, 
        f'(u) \left|\, \partial _x u \, \right|^{q+1} 
        \, \right) \\
&\qquad \; \, + \frac{q}{q+1} \, 
          f''(u) \left|\, \partial _x u \, \right|^{q+1}\partial _x u 
+\mu \, p\, q \left|\, \partial _x u \, \right|^{p+q-2} 
   \bigl( \, \partial _x^2 u \, \bigr)^2 
= 0. 
\end{aligned}
\end{align}
Integrating the divergence form (5.2) with respect to $x$, 
we have 
\begin{align}
\begin{aligned}
&\frac{1}{q+1} \, 
\frac{\mathrm{d}}{\mathrm{d}t}\, 
\,\Vert \, \partial _x u 
(t) \, \Vert_{L^{q+1}}^{q+1} 
+\mu \, p\, q \int ^{\infty }_{-\infty } 
\left|\, \partial _x u \, \right|^{p+q-2} 
\bigl( \, \partial _x^2 u \, \bigr)^2 
\, \mathrm{d}x \\
&\qquad \qquad \qquad \qquad \quad \; \; \; \: 
 +\frac{q}{q+1} \int ^{\infty }_{-\infty } 
  f''(u) \left|\, \partial _x u \, \right|^{q+1}\partial _x u 
  \, \mathrm{d}x 
= 0. 
\end{aligned}
\end{align}
Now we separate the integral region 
to the third term on the left-hand side of (5.3) as 
\begin{align}
\begin{aligned}
&\int ^{\infty }_{-\infty } 
 f''(u) \left|\, \partial _x u \, \right|^{q+1}\partial _x u 
 \, \mathrm{d}x \\
&= \int _{\partial _x u \geq 0 } + \int _{\partial _x u < 0 } \\
&= \int _{\partial _x u \geq 0 } 
   f''(u) \left|\, \partial _x u \, \right|^{q+2} \, \mathrm{d}x 
  - \int _{\partial _x u < 0 } 
  f''(u) \left|\, \partial _x u \, \right|^{q+2} \, \mathrm{d}x. 
\end{aligned}
\end{align}
Substituting (5.4) into (5.3), we get the following equality 
\begin{align}
\begin{aligned}
&\frac{1}{q+1} \, 
\frac{\mathrm{d}}{\mathrm{d}t}\, 
\,\Vert \, \partial _x u 
(t) \, \Vert_{L^{q+1}}^{q+1} 
+\mu \, p\, q \int ^{\infty }_{-\infty } 
\left|\, \partial _x u \, \right|^{p+q-2} 
\bigl( \, \partial _x^2 u \, \bigr)^2 
\, \mathrm{d}x \\
& \; \: \, +\frac{q}{q+1} \int _{\partial _x u \geq 0 } 
   f''(u) \left|\, \partial _x u \, \right|^{q+2} \, \mathrm{d}x 
= \frac{q}{q+1} \int _{\partial _x u < 0 } 
  f''(u) \left|\, \partial _x u \, \right|^{q+2} \, \mathrm{d}x. 
\end{aligned}
\end{align}

We have the following result 
which 
plays the most important role in the proof of (5.1). 

\medskip

\noindent
{\bf Lemma 5.1.}\quad {\it
For each $q>1$, there exists a positive constant $C_{q}$ 
such that 
\begin{align}
\begin{aligned}
\int _{\partial _x u < 0 } 
f''(u) \left|\, \partial _x u \, \right|^{q+2} \, \mathrm{d}x 
\leq C_{q}
     \int _{\partial _x u < 0 } 
     \left|\, \partial _x \phi \, \right|^{q+2} \, \mathrm{d}x. 
\end{aligned}
\end{align}
}

\noindent
In fact, taking care of the relation 
\begin{align}
\begin{aligned}
\partial _x u = \partial _x \tilde{U} + \partial _x \phi <0 \, 
\Longleftrightarrow  \, \partial _x \phi <0, 
                     \, \partial _x \tilde{U} < \bigl|\, \partial _x \phi \, \bigr|, 
\end{aligned}
\end{align}
we immediately have 
\begin{align}
\begin{aligned}
\int _{\partial _x u < 0 } 
&f''(u) \left|\, \partial _x u \, \right|^{q+2} \, \mathrm{d}x \\
&\leq 2^{q+2} \, 
     \left( \, \sup_{0 \leq u \leq \widetilde{C} + 1} f''(u) \, \right)
     \int _{\partial _x \phi < 0, \partial _x \tilde{U}< | \partial _x \phi | } 
     \left|\, \partial _x \phi \, \right|^{q+2} \, \mathrm{d}x.
\end{aligned}
\end{align}

\medskip

\noindent
{\bf Remark 5.1.}\quad
Under the relation (5.7), 
noting 
$$
\int ^{\infty }_{0} \int ^{\infty }_{-\infty } 
\left|\, \partial _x \phi \, \right|^{p+1} \, \mathrm{d}x\mathrm{d}t < \infty
$$
from (4.1) in Section 4 
and taking $q=p-1$ to (5.5), 
we 
can easily show that 
for $p \geq \frac{3}{2}$, 
\begin{align}
\begin{aligned}
\frac{1}{p} \, 
\,\Vert \, \partial _x u 
(t) \, \Vert_{L^{p}}^{p} 
&+\mu \, p\, (p-1) \int ^{\infty }_{0} \int ^{\infty }_{-\infty } 
\left|\, \partial _x u \, \right|^{2p-3} 
\bigl( \, \partial _x^2 u \, \bigr)^2 
\, \mathrm{d}x\mathrm{d}t  \\
&+\frac{p-1}{p} \int ^{\infty }_{0} \int _{\partial _x u \geq 0} 
  f''(u) \left|\, \partial _x u \, \right|^{p+1} 
  \, \mathrm{d}x\mathrm{d}t  
< \infty , 
\end{aligned}
\end{align}
which namely means that for $p \geq \frac{3}{2}$, 
\begin{eqnarray}
\left\{\begin{array}{ll}
\displaystyle{ 
\frac{\mathrm{d}}{\mathrm{d}t}\, 
\,\Vert \, \partial _x u  \, \Vert_{L^{p}}^{p} } 
\in L_{t}^1(0,\infty), \\[10pt]
\displaystyle{ 
\int ^{\infty }_{-\infty } 
\left|\, \partial _x u \, \right|^{2p-3} 
\bigl( \, \partial _x^2 u \, \bigr)^2 
\, \mathrm{d}x } 
\in L_{t}^1(0,\infty),  \\[10pt]
\displaystyle{ 
\int ^{\infty }_{-\infty } 
f''(u) \left|\, \partial _x u \, \right|^{p+1} 
\, \mathrm{d}x } 
\sim 
\displaystyle{ 
\int _{u>0 } 
\left|\, \partial _x u \, \right|^{p+1} 
\, \mathrm{d}x } 
\in L_{t}^1(0,\infty). 
\end{array} 
\right.\,
\end{eqnarray}

\bigskip

Integrating (5.5) with respect to $t$ and taking $q=p$, 
we have the 
energy equality 
\begin{align}
\begin{aligned}
&\frac{1}{p+1} \, 
\,\Vert \, \partial _x u 
(t) \, \Vert_{L^{p+1}}^{p+1} 
+\mu \, p^2 \int ^{t }_{0 } \int ^{\infty }_{-\infty } 
\left|\, \partial _x u \, \right|^{2(p-1)} 
\bigl( \, \partial _x^2 u \, \bigr)^2 
\, \mathrm{d}x \mathrm{d}\tau \\
&\; \: \, 
 +\frac{p}{p+1} \int ^{t }_{0 } \int _{\partial _x u \geq 0 } 
   f''(u) \left|\, \partial _x u \, \right|^{p+2} \, \mathrm{d}x \mathrm{d}\tau \\
&= \frac{1}{p+1} \, 
   \,\Vert \, \partial _x u_{0} \, \Vert_{L^{p+1}}^{p+1} 
   + \frac{p}{p+1} \int ^{t }_{0 } \int _{\partial _x u < 0 } 
     f''(u) \left|\, \partial _x u \, \right|^{p+2} \, \mathrm{d}x \mathrm{d}\tau. 
\end{aligned}
\end{align}

\noindent
The 
most difficult term to stimate is the second term on the right-hand side. 
We 
prepare the following ``boundary zero condition type'' interpolation inequality 
to overcome the difficulty. 

\medskip

\noindent
{\bf Lemma 5.2.}\quad {\it
It holds that 
\begin{align}
\begin{aligned}
& \int _{\partial _x u < 0 } 
  \left|\, \partial _x u \, \right|^{p+2} \, \mathrm{d}x \\
& \le C_{p} 
     \left( \, 
     \int _{\partial _x u < 0 } 
     \left|\, \partial _x u \, \right|^{2(p-1)} 
     \bigl( \, \partial _x^2 u \, \bigr)^2 
     \, \mathrm{d}x \, 
     \right)^{\frac{1}{3p+1}} 
     \left( \, 
     \int _{\partial _x u < 0 } 
     \left|\, \partial _x u \, \right|^{p+1} 
     \, \mathrm{d}x \, 
     \right)^{\frac{3p+2}{3p+1}}. 
\end{aligned}
\end{align}
}

\medskip

\noindent
{\bf Proof of Lemma 5.2.}\quad 
Since $ \partial _x u $ is absolutely continuous, 
we first note that 
for any 
$
x \in \bigl\{ \, x \in \mathbb{R} \, \, \bigr. 
 \bigl| \, \partial _x u \, < \, 0 \, \bigr\}, 
$
there exsists 
$
x_{k} \in \mathbb{R}\cup \{ - \infty \}
$
such that 
$$
\partial _x u (x_{k}) = 0, \; \; 
  \partial _x u (y) \, < \, 0 \; 
  \bigl( \, y \in 
  ( x_{k},x ) \, \bigr).
$$
Therefore by using the Cauchy-Schwarz inequality, 
it follows that for such $x$ and $x_{k}$ with $q\geq p\, (\, >1\, )$, 
\begin{align}
\begin{aligned}
&\left|\, \partial _x u \, \right|^q 
= \left(\, - \partial _x u \, \right)^q 
= q \int_{x_k}^{x} 
   \left(\, - \partial _x u \, \right)^{q-1} \left(\, - \partial _x^2 u \, \right) 
   \, \mathrm{d}y \\
&\leq q\, 
      \int _{ \partial _x u < 0 } 
      \left(\, - \partial _x u \, \right)^{q-1} \left(\, - \partial _x^2 u \, \right) 
      \, \mathrm{d}x \\
&\leq q\, \left( \,  
      \int _{ \partial _x u < 0 } 
      \left(\, - \partial _x u \, \right)^{2(p-1)} 
      \left(\, - \partial _x^2 u \, \right)^2 
      \, \mathrm{d}x \, 
      \right)^{\frac{1}{2}} 
      \left( \, 
      \int _{ \partial _x u < 0 }   
      \left(\, - \partial _x u \, \right)^{2(q-p)} 
      \, \mathrm{d}x \, 
      \right)^{\frac{1}{2}}. 
\end{aligned}
\end{align}
Hence 
\begin{align}
\begin{aligned}
&\Vert \, \partial _x u (t
) \, \Vert
_{
 L_{x}^
 {\infty}\left( \left\{ \partial _x u < 0 \right\} \right)
 } \\
&\leq q^{\frac{1}{q}}\, 
      \left( \,  
      \int _{ \partial _x u < 0 } 
      \left|\, \partial _x u \, \right|^{2(p-1)} 
      \bigl( \, \partial _x^2 u \, \bigr)^2 
      \, \mathrm{d}x \, 
      \right)^{\frac{1}{2q}} 
      \left( \, 
      \int _{ \partial _x u < 0 }  
      \left|\, \partial _x u \, \right|^{2(q-p)} 
      \, \mathrm{d}x \, 
      \right)^{\frac{1}{2q}}. 
\end{aligned}
\end{align}
So we get 
\begin{align}
\begin{aligned}
&\int _{ \partial _x u < 0 }  
 \left|\, \partial _x u \, \right|^{p+2} 
 \, \mathrm{d}x \\
&\leq \Vert \, \partial _x u \, \Vert
      _{
      L_{x}^
      {\infty}\left( \left\{ \partial _x u < 0 \right\} \right)
      }
      \int _{ \partial _x u < 0 }   
      \left|\, \partial _x u \, \right|^{p+1} 
      \, \mathrm{d}x \\
&\leq q^{\frac{1}{q}}\, 
      \left( \,  
      \int _{ \partial _x u < 0 }  
      \left|\, \partial _x u \, \right|^{2(p-1)} 
      \bigl( \, \partial _x^2 u \, \bigr)^2 
      \, \mathrm{d}x \, 
      \right)^{\frac{1}{2q}} \\
& \qquad \qquad \; \: \, \,  \times 
      \left( \, 
      \int _{ \partial _x u < 0 } 
      \left|\, \partial _x u \, \right|^{2(q-p)} 
      \, \mathrm{d}x \, 
      \right)^{\frac{1}{2q}} 
      \left(\, 
      \int _{ \partial _x u < 0 }   
      \left|\, \partial _x u \, \right|^{p+1} 
      \, \mathrm{d}x \, 
      \right). 
\end{aligned}
\end{align}
Taking $q=\frac{3p+2}{2}$ in (5.15), 
we have 
\begin{align}
\begin{aligned}
& \int _{ \partial _x u < 0 }  
 \left|\, \partial _x u \, \right|^{p+2} \, \mathrm{d}x 
\le \left( \, \frac{3p+2}{2} \, \right)^{\frac{2}{3p+2}} 
    \left( \, 
    \int _{ \partial _x u < 0 }   
    \left|\, \partial _x u \, \right|^{p+1} 
    \, \mathrm{d}x \, 
    \right)^{\frac{3p+2}{3p+1}} \\
& \qquad \qquad \qquad \qquad \quad \; \; \: \: \, \,  \times 
  \left( \, 
  \int _{ \partial _x u < 0 }  
  \left|\, \partial _x u \, \right|^{2(p-1)} 
  \bigl( \, \partial _x^2 u \, \bigr)^2 
  \, \mathrm{d}x \, 
  \right)^{\frac{1}{3p+1}}. 
\end{aligned}
\end{align}
Thus we complete the proof. 

\bigskip

Using Young's inequality to Lemma 5.2, (5.12), we also have 

\medskip

\noindent
{\bf Lemma 5.3.}\quad{\it
It follows that for any $\epsilon>0$, 
there exists a positive constant $C_p({\epsilon})$ such that, 
\begin{align}
\begin{aligned}
& \int _{\partial _x u < 0 } 
  \left|\, \partial _x u \, \right|^{p+2} \, \mathrm{d}x \\
& \le \epsilon 
      \int _{\partial _x u < 0 } 
      \left|\, \partial _x u \, \right|^{2(p-1)} 
      \bigl( \, \partial _x^2 u \, \bigr)^2 
      \, \mathrm{d}x \, 
      + C_p({\epsilon}) \left( \, 
        \int _{\partial _x u < 0 } 
        \left|\, \partial _x u \, \right|^{p+1} 
        \, \mathrm{d}x \, 
        \right)^{\frac{3p+2}{3p}}. 
\end{aligned}
\end{align}
}
\medskip

\noindent
Substituting (5.17) with 
$\epsilon = \frac{\mu \, p^2}{2} $
into (5.11), we have 
\begin{align}
\begin{aligned}
&\frac{1}{p+1} \, 
\,\Vert \, \partial _x u 
(t) \, \Vert_{L^{p+1}}^{p+1} 
+\frac{\mu \, p^2}{2} \int ^{t }_{0 } \int ^{\infty }_{-\infty } 
\left|\, \partial _x u \, \right|^{2(p-1)} 
\bigl( \, \partial _x^2 u \, \bigr)^2 
\, \mathrm{d}x \mathrm{d}\tau \\
&\; \: \, 
 +\frac{p}{p+1} \int ^{t }_{0 } \int _{\partial _x u \geq 0 } 
   f''(u) \left|\, \partial _x u \, \right|^{p+2} \, \mathrm{d}x \mathrm{d}\tau \\
&\leq \frac{1}{p+1} \, 
   \,\Vert \, \partial _x u_{0 } 
   \, \Vert_{L^{p+1}}^{p+1} 
+ C_p 
  \int ^{t }_{0 } 
  \left( \, 
  \int _{\partial _x u < 0 } 
  \left|\, \partial _x u \, \right|^{p+1} 
  \, \mathrm{d}x \, 
  \right)^{\frac{2}{3p}+1} \mathrm{d}\tau . 
\end{aligned}
\end{align}

\noindent
Now recalling Lemma 5.1, we have 
\begin{align}
\begin{aligned}
\int _{\partial _x u < 0 } 
\left|\, \partial _x u \, \right|^{p+1} 
\, \mathrm{d}x 
&\leq C_{p} 
      \int ^{\infty}_{-\infty} 
      \bigl|\, \partial _x \phi \, \bigr|^{p+1} 
      \, \mathrm{d}x 
\in L_t^1(0,\infty). 
\end{aligned}
\end{align}
We also note $\frac{2}{3p}<1$ and focus on the fact
\begin{align}
\begin{aligned}
\left( \, 
\int _{\partial _x u < 0 } 
\left|\, \partial _x u \, \right|^{p+1} 
\, \mathrm{d}x \, 
\right)^{\frac{2}{3p}} 
\leq C_p \left( \, 1 +
     \int _{\partial _x u < 0 } 
     \left|\, \partial _x u \, \right|^{p+1} 
     \, \mathrm{d}x \, 
     \right)
\end{aligned}
\end{align}
for some posisive constant $C_p$. 
Hence, substituting (5.19) and (5.20) into (5.18), 
we have 
\begin{align}
\begin{aligned}
&\frac{1}{p+1} \, 
\,\Vert \, \partial _x u 
(t) \, \Vert_{L^{p+1}}^{p+1} 
+\frac{\mu \, p^2}{2} \int ^{t }_{0 } \int ^{\infty }_{-\infty } 
\left|\, \partial _x u \, \right|^{2(p-1)} 
\bigl( \, \partial _x^2 u \, \bigr)^2 
\, \mathrm{d}x \mathrm{d}\tau \\
&\; \: \, 
 +\frac{p}{p+1} \int ^{t }_{0 } \int _{\partial _x u \geq 0 } 
   f''(u) \left|\, \partial _x u \, \right|^{p+2} \, \mathrm{d}x \mathrm{d}\tau \\
&\leq \frac{1}{p+1} \, 
   \,\Vert \, \partial _x u_{0} \, \Vert_{L^{p+1}}^{p+1} 
+ C_p 
  \int ^{t }_{0 } 
  \int ^{\infty}_{-\infty} 
  \bigl|\, \partial _x \phi \, \bigr|^{p+1} 
  \, \mathrm{d}x 
  \mathrm{d}\tau \\
&\qquad \qquad \qquad \quad \quad \; \; \: \: \: \, 
            + C_p 
            \int ^{t }_{0 } 
            \,\Vert \, \partial _x u 
            (\tau) \, \Vert_{L^{p+1}}^{p+1} 
            \left( \, 
            \int ^{\infty}_{-\infty} 
            \bigl|\, \partial _x \phi \, \bigr|^{p+1} 
            \, \mathrm{d}x \, 
            \right) 
            \mathrm{d}\tau. 
\end{aligned}
\end{align}
By using 
Gronwall's inequality, we have 
\begin{align}
\begin{aligned}
\,\Vert \, \partial _x u 
(t) \, \Vert_{L^{p+1}}^{p+1} 
&\leq C_p \bigl( \, 
      \| \, \partial _x u_0 \, \|_{L^{p+1}}^{p+1} + \| \, \phi_0 \, \|_{L^2}^2 + 1 
      \, \bigr) \\
&\quad \times \exp 
 \left\{C_p  
 \int ^{\infty }_{0 } 
 \int ^{\infty}_{-\infty} 
 \bigl|\, \partial _x \phi \, \bigr|^{p+1} 
 \, \mathrm{d}x 
 \mathrm{d}t
 \right\}
< \infty. 
\end{aligned}
\end{align}
Hence, substituting (5.22) into (5.21), it finally holds 
\begin{align}
\begin{aligned}
&\| \, \partial _x u 
(t) \, \| _{L^{p+1}}^{p+1} 
+\int _0^t \int _{-\infty}^{\infty} 
 \bigl| \, \partial _x u \, \bigr|^{2(p-1)} 
 \left( \, \partial _x^2 u \, \right) 
 \, \mathrm{d}x \mathrm{d}\tau \\
&\qquad \qquad \quad \; \: \: \, 
 +\int _0^t 
  \| \, 
  \partial _x u 
  (\tau)
  \, \| _{L^{p+2}\left( \left\{x \in \mathbb{R}\, |\, u>0 \right\} \right)}^{p+2} 
  \, \mathrm{d}\tau \\
&\leq 
       C\bigl( \, \| \, \phi_0 \, \|_{L^2}, \: 
       \| \, \partial _x u_0 \, \|_{L^{p+1}} 
       \, \bigr). 
\end{aligned}
\end{align}
Thus, we do complete the proof of Proposition 3.2.

\bigskip 

\noindent
\section{Asymptotic behavior}
In this section, we shall obtain the asymptotic behavior 
\begin{equation}
\displaystyle{
\sup_{x \in \mathbb{R}}\left|\,  \phi(t,x) \, \right|
\xrightarrow[\, \, t\to  \infty \, ]
\ \ 0} 
\end{equation}
by utilizing Proposition 3.1, (3.11) and Proposition 3.2, (3.12). 
Noting the Gagliardo-Nirenberg inequality 
(cf. \cite{gag}, \cite{nash}, \cite{nir}), we have 
\begin{align}
\begin{aligned}
&\displaystyle{
 \sup_{x \in \mathbb{R}}\left|\,  \phi(t,x) \, \right|} \\
&\leq C_p\, 
      \| \, \phi 
      (t)\, \| _{L^{2}}^{\frac{2p}{3p+1}} 
      \bigl|\bigl| 
      \, \partial _x \phi 
      (t)\, 
      \bigr|\bigr| _{L^{p+1}}^{\frac{p+1}{3p+1}} \\
&\leq C_p\, 
      \bigl( \, \| \, \phi_0 \, \|_{L^2}^2 + 1 \, \bigr)^{\frac{2p}{3p+1}} 
      \left( \, 
      \bigl|\bigl| 
      \, \partial _x u 
      (t)\, 
      \bigr|\bigr| _{L^{p+1}}^{\frac{p+1}{3p+1}}
      + \bigl|\bigl| 
        \, \partial _x \tilde{U} 
      (t) \, 
        \bigr|\bigr| _{L^{p+1}}^{\frac{p+1}{3p+1}} \, 
      \right), 
\end{aligned}
\end{align}
and we also note 
\begin{equation}
\bigl|\bigl| 
\, \partial _x \tilde{U} 
(t) \, 
\bigr|\bigr| _{L^{p+1}}
\leq C_p\, (1+t)^{-\frac{1}{p+1}} 
\xrightarrow[\, \, t\to  \infty \, ]
\ \ 0.
\end{equation}
Hence, it suffices to prove 

\medskip

\noindent
{\bf Lemma 6.1.}\quad{\it 
It holds that 
\begin{equation}
\displaystyle{
\bigl|\bigl| 
\, \partial _x u 
(t) \, 
\bigr|\bigr| _{L^{p+1}}
\xrightarrow[\, \, t\to  \infty \, ]
\ \ 0}. 
\end{equation}
}

\medskip

\noindent
To show Lemma 6.1, we claim 

\medskip

\noindent
{\bf Lemma 6.2.}\quad{\it 
It holds that 
\begin{equation}
\frac{\mathrm{d}}{\mathrm{d}t}\, 
\bigl|\bigl| 
\, \partial _x u  \, 
\bigr|\bigr| _{L^{p+1}}^{p+1} 
\in L_t^1(0, \infty). 
\end{equation}
}

\medskip

\noindent
Once 
it holds, Lemma 6.1 immediately follows. 
In fact, for any sequence
$$
\left\{ \, t_k \, \right\}_{k=1}^{\infty} \subset [\, 0, \infty) \; \; 
\bigl( \, 
t_k \nearrow \infty \; (k\rightarrow \infty) \, 
\bigr), 
$$
we have 
\begin{align}
\begin{aligned}
&\left| \, 
 \bigl|\bigl| 
 \, \partial _x u 
 (t_m) \, 
 \bigr|\bigr| _{L^{p+1}}^{p+1} 
 - \bigl|\bigl| 
   \, \partial _x u 
   (t_n) \, 
   \bigr|\bigr| _{L^{p+1}}^{p+1} \, 
 \right| \\ 
&\leq \Biggl| \, 
      \int_{0}^{t_m} 
      \left| \, 
      \frac{\mathrm{d}}{\mathrm{d}t}\, 
      \bigl|\bigl| 
      \, \partial _x u 
      (t) \, 
      \bigr|\bigr| _{L^{p+1}}^{p+1} \, 
      \right| \, 
      \mathrm{d}t 
      - \int_{0}^{t_n} 
        \left| \, 
        \frac{\mathrm{d}}{\mathrm{d}t}\, 
        \bigl|\bigl| 
        \, \partial _x u 
        (t) \, 
        \bigr|\bigr| _{L^{p+1}}^{p+1} \, 
        \right| \, 
        \mathrm{d}t
     \, 
     \Biggr| \\ 
& \xrightarrow[\, \, m, \: n \to  \infty \, ]
  \ \ 0. 
\end{aligned}
\end{align}
Therefore 
$
\bigl\{ \, 
\Vert 
\, \partial _x u 
(t_k) \, 
\Vert _{L^{p+1}}^{p+1} \, 
\bigr\}_{k=1}^{\infty} 
\subset \mathbb{R}
$
is a Cauchy sequence 
which has a limit $\alpha$ in $\mathbb{R}$. 
Because 
$\left\{ \, t_k \, \right\}_{k=1}^{\infty}$ is arbitrarily taken, 
the limit $\alpha$ should be independent of $\left\{ \, t_k \, \right\}_{k=1}^{\infty}$ 
satisfying 
\begin{equation}
\alpha = \displaystyle{
         \lim _{t \rightarrow \infty }
         \bigl|\bigl| 
         \, \partial _x u 
         (t) \, 
         \bigr|\bigr| _{L^{p+1}}^{p+1}
         }
       = \displaystyle{
            \lim _{t \rightarrow \infty }
            \bigl|\bigl| 
            \, \partial _x \phi 
            (t)  \, 
            \bigr|\bigr| _{L^{p+1}}^{p+1}
            }. 
\end{equation}
Thus, it holds $\alpha =0$ by 
$
\Vert  
\, \partial _x \phi \, 
\Vert  _{L^{p+1}}^{p+1}
\in L_t^1(0,\infty).
$

It remains to prove Lemma 6.2. 

\medskip

\noindent
{\bf Proof of Lemma 6.2.}\quad 
Direct calculation shows that 
\begin{align}
\begin{aligned}
&\frac{\mathrm{d}}{\mathrm{d}t}\, 
 \int _{-\infty}^{\infty} 
 \bigl|
 \, \partial _x u  \, 
 \bigr|^{p+1} \, 
 \mathrm{d}x \\ 
&= (p+1) \int _{-\infty}^{\infty} 
   \Bigl( \, \bigl|
   \, \partial _x u  \, 
   \bigr|^{p-1} \, 
   \partial _x u \, \Bigr) \, 
   \partial _x \bigl( \, \partial _t u \, \bigr) \, 
   \mathrm{d}x \\
&= - \int _{-\infty}^{\infty} 
   \left( \, 
   p\, f''(u)\, \bigl|\, \partial _x u  \, \bigr|^{p+1} \, 
   \partial _x u 
   + \mu \, p^2 (p+1) 
     \bigl|\, \partial _x u  \, \bigr|^{2(p-1)} 
     \bigl( \, \partial _x^2 u \, \bigr)^2 \, 
   \right) \, 
   \mathrm{d}x. 
\end{aligned}
\end{align}
Integrating (6.8) with respect to $t$, we have by (3.11) and (3.12), 
\begin{align}
\begin{aligned}
&\int ^{\infty }_{0 } 
 \left| \, 
 \frac{\mathrm{d}}{\mathrm{d}t}\, 
 \int _{-\infty}^{\infty} 
 \bigl|
 \, \partial _x u  \, 
 \bigr|^{p+1} \, 
 \mathrm{d}x \, 
 \right| \, 
 \mathrm{d}t \\
& \leq p \, 
       \left( \, \sup_{0 \leq u \leq \widetilde{C} + 1} f''(u) \, \right) 
       \int ^{\infty }_{0 } \int _{u>0}
       \bigl|\, \partial _x u  \, \bigr|^{p+2} \, 
       \mathrm{d}x \mathrm{d}t \\ 
& \qquad + \mu \, p^2 (p+1) 
           \int ^{\infty }_{0 } \int _{-\infty}^{\infty} 
           \bigl|\, \partial _x u  \, \bigr|^{2(p-1)} 
           \bigl( \, \partial _x^2 u \, \bigr)^2 \, 
           \mathrm{d}x \mathrm{d}t \\ 
& \leq 
        C_{p}\bigl( \, 
        \| \, \phi_0 \, \|_{L^2}, \: \| \, \partial _x u_0 \, \|_{L^{p+1}} 
        \, \bigr). 
\end{aligned}
\end{align}
Hence, we complete the proof of Lemma 6.1. 

Thus we have obtained the asymptotic behavior (6.1).

\bigskip

{\bf Acknowledgement.}\quad
The auther 
thanks 
Professor Akitaka Matsumura for 
his constant help and guidance. 

\end{document}